\documentclass{article}
\usepackage{graphicx}
\usepackage{amsthm}
\usepackage{amsmath}
\usepackage{amsfonts}
\usepackage{amssymb}

\newtheorem{theorem}{Theorem}

\newtheorem{corollary}{Corollary}

\newtheorem*{definition}{Definition}

\newtheorem*{lemma}{Lemma}
\newtheorem{notation}{Notation}
\newtheorem*{note}{Note}

\newtheorem{proposition}{Proposition}
\newtheorem{remark}{Remark}

\begin{document}

\title{On certain families of naturally graded Lie algebras\thanks{%
Research partially supported by the D.G.I.C.Y.T project PB98-0758}}
\author{Jos\'{e} Mar\'{\i}a Ancochea Berm\'{u}dez\thanks{%
corresponding author: e-mail: Jose\_Ancochea@mat.ucm.es} \and Rutwig
Campoamor Stursberg\thanks{e-mail: rutwig@nfssrv.mat.ucm.es} \\
%EndAName
Departamento de Geometr\'{\i}a y Topolog\'{\i}a\\
Fac. CC. Matem\'{a}ticas Univ. Complutense\\
28040 Madrid ( Spain )}
\date{}
\maketitle

\begin{abstract}
In this work large families of naturally graded nilpotent Lie algebras in arbitrary 
dimension and characteristic sequence $\left(n,q,1\right)$ with $n\equiv 1
\left(mod 2\right)$ satisfying the centralizer property are given. This
centralizer property constitutes a generalization, for any nilpotent algebra,
of the structural properties characterizing the Lie algebra $Q_{n}$. By
considering certain cohomological classes of the space
$H^{2}\left(\frak{g},\mathbb{C}\right)$, it is shown that, with few
exceptions, the isomorphism classes of these algebras are given by central
extensions of $Q_{n}$ by $\mathbb{C}^{p}$ which preserve the nilindex and the
natural graduation.\newline Math. Subj. Class : 17B30, 17B56, 17B70.\newline
\textit{keywords : nilpotent, naturally graded, extension, cohomology}.
\end{abstract}

\section*{Introduction}
 The first systematic results
about naturally graded nilpotent Lie were obtained in the sixties in the
context of the analysis of the variety of nilpotent Lie algebra laws. In
$1966$ Vergne [9] concentrated on the less nilpotent Lie algebras, which she
called filiform, and classified the naturally graded Lie algebras having this
property. This gave the key to a first estimation of the number of irreducible
components of the variety [10]. The classification result was, in a certain
manner, surprising : there are only two models, called respectively $L_{n}$
and $Q_{n}$, where the second exists only in even dimension. The filiform
model Lie algebra, as $L_{n}$ is called usually, is without doubt the most and
best studied nilpotent Lie algebra over the last thirty years. Most studies
dedicated to filiform Lie algebras and its deformations are dedicated to this
model. For this reason, the algebra $Q_{n}$ has often been relegated to a
secondary position. However, the nonexistence of this algebra in odd
dimensions makes it, from a structural point of view, much more interesting
than the model algebra $L_{n}$. In fact, when considering the centralizers
$C_{Q_{n}}I$ of ideals $I$ of the central descending sequence in $Q_{n}$, we
find an intriguing property, namely, that the integer part of the nilindex
modulus $2$ works as frontier between those ideals of the sequence contained
in its associated centralizer in $Q_{n}$ and those whose intersection of the
centralizer with the complementary of the ideal in the algebra is nonzero.
This fact can be used to estimate how far is $Q_{n}$ from having an abelian
commutator algebra; the index above tells that this algebra is as far as
posiible from having it. This question leads naturally to search for a
generalization of this property for other naturally graded nilpotent Lie
algebras. Making use of the characteristic sequence of an algebra, we can
concentrate on concrete classes of algebras. The cohomology space
$H^{2}\left(\frak{g},\mathbb{C}\right)$ of $\frak{g}$ with values in the base
field are of wide interest to determine particular classes of central
extensions of $\frak{g}$ by $\mathbb{C}$ which preserve either the natural
graduation or any other property of the extended algebra. By introducing a
partition of this space, we are in situation of isolating the cohomology
classes which make direct reference to the property of the centralizers. This
can be used to achieve a complete classification, for a fixed characteristic
sequence, of algebras bahaving as $Q_{n}$ does. This is done for the sequences
$\left(2m-1,1,1\right)$ and $\left(2m-1,2,1\right)$ in arbitrary dimension,
and can be applied to any other sequence. Our main purpose is, however, to
provide families, in arbitrary dimensions and characteristic sequences
$\left(n,q,1\right)$ with $n\equiv 1 \left(mod 2\right)$ and $q\geq 1$, of
naturally graded nilpotent Lie algebras with this centralizer property. The
obtained algebras can be interpreted as the analogue, for its corresponding
characteristic sequence, of the Lie algebra $Q_{n}$.

\section{Preliminaries and notations}
 
 Whenever we speak about Lie algebras
in this work, we refer to finite dimensional complex Lie algebras.\bigskip
 
\begin{definition}
 Let $\frak{g}$ be a finite dimensional vectorial space
over $\mathbb{C}$.
 A Lie algebra law over $\mathbb{C}^{n}$ is a bilinear
alternated mapping $\mu\in Hom\left( 
\mathbb{C}^{n}\times\mathbb{C}^{n},\mathbb{C}^{n}\right)$ which satisfies the
conditions
 
 \begin{enumerate}
 \item $\mu\left(  X,X\right) 
=0,\;\forall\;X\in\mathbb{C}^{n}$
 
 \item $\mu\left(  X,\mu\left(  Y,Z\right)
 \right)  +\mu\left(  Z,\mu\left(
 X,Y\right)  \right)  +\mu\left( 
Y,\mu\left(  Z,X\right)  \right)
 =0,\;\forall\;X,Y,Z\in\mathbb{C}^{n}$, \
\newline  ( Jacobi identity )
 \end{enumerate}
 If  $\mu$ is a Lie algebra
law, the pair  $\frak{g}= \left(\mathbb{C}^{n}, \mu \right)$ is  called Lie
algebra. From now on we identify the Lie algebra with its law $\mu$. 
\end{definition}
 
 \begin{remark}
 We say that $\mu$ is the law of
$\frak{g}$, and where necessary we use the bracket notation to describe the
law :
 \[
 \left[  X,Y\right]  =\mu\left(  X,Y\right) 
,\;\forall\;X,Y\in\frak{g}%
 \]
 The nondefined brackets are zero or obtained
by antisymmetry. 
 \end{remark}
 
 \bigskip To any Lie algebra we can
associate the following sequence :
 
 \begin{align*}
 C^{0}\frak{g}  & 
=\frak{g}\supset C^{1}\frak{g}=\supset
 C^{2}\frak{g}=\left[ 
C^{1}\frak{g},\frak{g}\right]  \supset...\supset
 C^{k}\frak{g}=\left[ 
C^{k-1}\frak{g},\frak{g}\right]  \supset...
 \end{align*}
 called the
descending central sequence of $\frak{g}$.
 
 \begin{definition}
 A Lie
algebra $\frak{g}$ is called nilpotent if there exists an integer ( called
nilindex $n\left( \frak{g}\right)$ of $\frak{g}$) $k\geq1$ such that 
$C^{k}\frak{g}=\left\{
 0\right\}$ and $C^{k-1}\frak{g}=\left\{0\right\}$.
\end{definition}
 
 \begin{definition}
 An $n$-dimensional nilpotent Lie
algebra is called filiform if 
 \[
 \dim C^{k}\frak{g}=n-k-1,\;1\leq k\leq
n-1
 \]
 \end{definition}
 
 \begin{remark}
 Calling
$p_{i}=dim\left(\frac{C^{i-1}\frak{g}}{C^{i}\frak{g}}\right)$ for $1\leq\leq
n\left( \frak{g}\right)$, the type of the nilpotent Lie algebra is the
sequence $\left\{p_{1},..,p_{r}\right\}$. Then a filiform algebra corresponds
to those of type $\left\{2,1,..,1\right\}$ [10].
 \end{remark}
 
 We recall
the laws for the $\left(n+1\right)$-dimensional filiform Lie algebras $L_{n}$
and $Q_{n}$, which are the only ones we will use here :
 
 \begin{enumerate}
\item $L_{n}$ $\left(  n\geq3\right)  :$
 \[
 \lbrack
X_{1},X_{i}]=X_{i+1},\;2\leq i\leq n
 \]
 over the basis
$\left\{X_{1},..,X_{n+1}\right\}$.
 
 \item $Q_{2m-1}$ $\left(  m\geq3\right) 
:$%
 \begin{align*}
 \lbrack X_{1},X_{i}]  & =X_{i+1},\;1\leq i\leq2m-1\\
\lbrack X_{j},X_{2m+1-j}]  & =\left(  -1\right)  ^{j}X_{2m},\;2\leq j\leq m
\end{align*}
 over the basis $\left\{X_{1},..,X_{2m}\right\}$.
\end{enumerate}

 \begin{definition}
 A Lie algebra $\frak{g}$ \ is graded
over  $\mathbb{Z}$ if
 it admits a decomposition 
 \[
\frak{g}=\bigoplus_{k\in\mathbb{Z}}\frak{g}_{k}%
 \]
 where the $\frak{g}_{k}$
are $\mathbb{C}$-subspaces of
 $\frak{g}$ which satisfy  $\left[ 
\frak{g}_{r},\frak{g}_{s}\right]
 \subset\frak{g}_{r+s}$ ,
$\;r,s\in\mathbb{Z}$.
 \end{definition}
 
 Observe that any graduation defines
a sequence
 \[
 S_{k}=F_{k}\left(  \frak{g}\right)  =\bigoplus_{t\geq
k}\frak{g}_{t}%
 \]
 with the properties
 
 \begin{enumerate}
 \item
$\frak{g}=%
 %TCIMACRO{\dbigsqcup }%
 %BeginExpansion
{\displaystyle\bigsqcup}
 %EndExpansion
 \,S_{k}$

\item $\left[  S_{i},S_{j}\right]  \subset S_{i+j}\;\forall i,j$

\item $S_{i}\subset S_{j}$ \ si $i>j$
\end{enumerate}

\begin{definition}
A family $\left\{  S_{i}\right\}  $ of subspaces of 
$\frak{g}$ define a filtration ( descending ) over $\frak{g}$ if
it satisfies properties $1), 2), 3)$. The algebra is called filtered.
\end{definition}

The construction can be reversed, i.e., any filtration defines a graduation by taking  $\frak{g}
_{k}=\frac{S_{k-1}}{S_{k}}$ \ for $k\geq1$. The graduation is called associated to the filtration  $\left\{  S_{i}\right\}  $ and it defines a Lie algebra with the rule
\[
\left[  \overset{-}{X_{i}},\overset{-}{X_{j}}\right]  =\overline{\left[
X_{i},X_{j}\right]  },\;
\]
where $\overset{-}{X_{k}}=X_{k}$ $\left(  \operatorname{mod}\;S_{k+1}\right)
$ for  $k=i,j$ \ and \  $\overline{\left[  X_{i},X_{j}\right]  }\,=\left[
X_{i},X_{j}\right]  \;\left(  \operatorname{mod}\;S_{i+j+1}\right)  $

\begin{definition}
A nilpotent Lie algebra is called naturally graded if 
$\frak{g}\simeq\frak{gr}\left(  \frak{g}\right)$, where $\frak{gr}\left(\frak{g}\right)$ is the graduation associated to the filtration induced in $\frak{g}$ by the central descending sequence.
\end{definition}

It follows immediately that both $L_{n}$ and $Q_{n}$ are naturally graded. They are in fact the only filiform Lie algebras having this property [10].\bigskip

Let $\frak{g}_{n}=\left(  \mathbb{C}^{n},\mu\right)  $ be a nilpotent Lie algebra. 
For any nonzero vector $X\in\frak{g}_{n}-C^{1}\frak{g}_{n}$
let $c\left(  X\right)  $ be the ordered sequence of a similitude invariant for the nilpotent operator $ad_{\mu}\left(  X\right)
\dot{,}$ i.e., the ordered sequence of dimensions of Jordan blocks of this operator. The set of these sequences is ordered lexicographically.

\begin{definition}
The characteristic sequence of  $\frak{g}_{n}$ is an isomorphism invariant  $c\left(  \frak{g}_{n}\right)  $ defined by 
\[
c\left(  \frak{g}_{n}\right)  =\max_{X\in\frak{g}_{n}-C^{1}\frak{g}_{n}%
}\left\{  c\left(  X\right)  \right\}
\]
A nonzero vector $X\in\frak{g}_{n}-C^{1}\frak{g}_{n}$  for which 
$c\left(  X\right)  =c\left(  \frak{g}_{n}\right)  $ 
is called characteristic vector.
\end{definition}

\begin{remark}
In particular, the algebras with maximal characteristic sequence 
$\left(  n-1,1\right)  $ correspond to the filiform algebras introduced by Vergne.
\end{remark}

It is often convenient to use the so called ontragradient representation of a Lie algebra $\frak{g}$. Let $n=dim\left(\frak{g}\right)$ and $\left\{X_{1},..,X_{n}\right\}$ be a basis. If $C_{i,j}^{k}$ are the structure constants of the algebra law $\mu$, we can define, over the dual basis $\left\{\omega_{1},..,\omega_{n}\right\}$, the differential
\[
d_{\mu}\omega_{i}\left(X_{j},X_{k}\right)=-C_{jk}^{i}
\]
Then the Lie algebra is rewritten as 
\[
d\omega_{i}=-C_{jk}^{i}\omega_{j}\wedge\omega_{k}\; 1\leq i,j,k\leq n,
\]
The Jacobi condition is equivalent to $d^{2}\omega_{i}=0$ for all $i$.

\begin{remark}
In what follows we will use the preceding form, up to the sign, to describe Lie algebras. It will be seen that certain structural properties are better seen by using this form. For example, this was the procedure to analyze the Lie algebra models [4].
\end{remark}

\subsection{The spaces $H^{2}\left(\frak{g},\mathbb{C}\right)$}

Recall that the space $H^{2}\left(  \frak{g},\mathbb{C}^{p}\right)  $
can be interpreted as the space of classes of $p$-dimensional central extensions of the Lie algebra $\frak{g}$. We recall the elementary facts :

Let $\frak{g}$ be an $n$-dimensional nilpotent Lie algebra with law $\mu_{0}$. 
A central extension of $\frak{g}$ by $\mathbb{C}^{p}$ is an exact sequence of Lie algebras
\[
0\longrightarrow\mathbb{C}^{p}\longrightarrow\overset{-}{\frak{g}%
}\longrightarrow\frak{g}\longrightarrow0\text{ }%
\]
such that $\mathbb{C}^{p}\subset Z\left(  \overset{-}{\frak{g}}\right)  .$ Let
$\alpha$ be a cocycle of the De Rham cohomology $Z^{2}\left(  \frak{g},\mathbb{C}%
^{p}\right)  .$ This gives the extension
\[
0\longrightarrow\mathbb{C}^{p}\longrightarrow\mathbb{C}^{p}\oplus\frak{g}%
\longrightarrow\frak{g}\longrightarrow0
\]
with associated law $\mu=\mu_{0}+\alpha$ defined by%
\[
\mu\left(  \left(  a,x\right)  ,\left(  b,y\right)  \right)  =\left(
\alpha\mu_{0}\left(  x,y\right)  ,\mu_{0}(x,y\right)  )
\]
In the following we are only interested in extensions of $\mathbb{C}$ by $\frak{g}$, i.e, extensions of degree one. It is well known that the space of $2$-cocycles
$Z^{2}\left(  \frak{g},\mathbb{C}\right)$ is identified with the space of linear forms
over $\bigwedge^{2}\frak{g}$ which are zero over the subspace $\Omega$ :%
\[
\Omega:=\left\langle \mu_{0}\left(  x,y\right)  \wedge z+\mu_{0}\left(
y,z\right)  \wedge x+\mu_{0}\left(  z,x\right)  \wedge y\right\rangle
_{\mathbb{C}}%
\]
The extension classes are defined modulus the coboundaries 
$B^{2}\left(  \frak{g},\mathbb{C}\right)  .$ This allows to identify the cohomology space
$H^{2}\left(  \frak{g},\mathbb{C}\right)  $ with the dual of the space 
$\frac{Ker\;\lambda}{\Omega},$ where $\lambda\in Hom\left(  \bigwedge
^{2}\frak{g},\frak{g}\right)  $ is defined as
\[
\lambda\left(  x\wedge y\right)  =\mu_{0}\left(  x,y\right)\; x,y\in\frak{g}
\]
In fact we have $H_{2}\left(\frak{g},\mathbb{C}\right)=\frac{Ker\lambda}{\Omega}$ for the $2$-homology space, and as $H^{2}\left(\frak{g},\mathbb{C}\right)=Hom_{\mathbb{C}}\left(H_{2}\left(\frak{g},\mathbb{C}\right),\mathbb{C}\right)$ the assertion follows.

\begin{notation}
Let $\varphi_{ij}\in H^{2}\left(  \frak{g},\mathbb{C}\right)  $ the cocycles
defined by 
\[
\varphi_{ij}\left(  X_{k},X_{l}\right)  =\delta_{ik}\delta_{jl}
\]
\end{notation}

Observe that a cocycle  $\varphi$ can be written as a linear combination of the preceding cocycles. We have  :

\begin{lemma}
$\sum a^{ij}\varphi_{ij}=0$ if and only if  $\sum a^{ij}\left(  X_{i}\wedge
X_{j}\right)\in\Omega$
\end{lemma}

Let $\frak{g}$ be an $n$-dimensional nilpotent Lie algebra. The subspace of central extensions is noted by  $E_{c,1}\left(  \frak{g}\right)$. It has been shown that this space is 
irreducible and constructible. However, for our purpose this space is too general. We only need certain cohomology classes of this space.

\begin{notation}
For $k\geq2$ let %
\[
H_{k}^{2,t}\left(  \frak{g},\mathbb{C}\right)  =\left\{  \varphi_{ij}\in
H^{2}\left(  \frak{g},\mathbb{C}\right)  \;|\;i+j=2t+1+k\right\}  ,\;1\leq t\leq
\left[  \frac{n-3}{2}\right],\
\]

\[
H_{k}^{2,\frac{t}{2}}\left(  \frak{g},\mathbb{C}\right)  =\left\{
\varphi_{ij}\in H^{2}\left(  \frak{g},\mathbb{C}\right)
\;|\;i+j=t+1+k\right\}  ,\;t\in\{1,..,\left[  \frac{n-3}{2}\right]
\},t\equiv1\left(  \operatorname{mod}\,2\right)
\]

\end{notation}

These cocycles are essential to determine the central extensions which are additionally naturally graded. If $\mathbf{E}_{c,1}\left(\frak{g}\right)  $ denotes the central extensions that are naturally graded, we consider the subspaces 
\[
\mathbf{E}_{c,1}^{t,k_{1},..,k_{r}}\left(  \frak{g}\right)  =\left\{  \mu
\in\mathbf{E}_{c,1}\left(  \frak{g}\right)  \;|\;\mu=\mu_{0}+\left(
%TCIMACRO{\tsum }%
%BeginExpansion
{\textstyle\sum}
%EndExpansion
\varphi_{ij}^{k_{i}}\right)  ,\;\varphi_{ij}^{k_{i}}\in H_{k_{i}}^{2,t}\left(
\frak{g},\mathbb{C}\right)  \right\}
\]%

\[
\mathbf{E}_{c,1}^{\frac{t}{2},k_{1},..,k_{r}}\left(  \frak{g}\right)
=\left\{  \mu\in\mathbf{E}_{c,1}\left(  \frak{g}\right)  \;|\;\mu=\mu
_{0}+\left(
%TCIMACRO{\tsum }%
%BeginExpansion
{\textstyle\sum}
%EndExpansion
\varphi_{ij}^{k_{i}}\right)  ,\;\varphi_{ij}^{k_{i}}\in H_{k_{i}}^{2,\frac{t}{2}%
}\left(  \frak{g},\mathbb{C}\right)  \right\}
\]
where $0\leq k_{j}\in\mathbb{Z},$ $j=1,..,r.$\newline
Given a basis $\left\{X_{1},..,X_{n},X_{n+1}\right\}$ of $\mu$ belonging to any of these spaces, the Lie algebra law is defined by : 

\[
\mu\left(  X_{i},X_{j}\right)  =\mu_{0}\left(  X_{i},X_{j}\right)  +\left(
\sum\varphi_{ij}^{k}\right)  X_{n+1},\;1\leq i,j\leq n\; 
\]

\begin{lemma}
As vector spaces, the following identity holds :
\[
\mathbf{E}_{c,1}\left(  \frak{g}\right)  =\sum_{t,k}\mathbf{E}_{c,1}%
^{t,k_{1},..,k_{r}}\left(  \frak{g}\right)  +\mathbf{E}_{c,1}^{\frac{t}%
{2},k_{1},..,k_{r}}\left(  \frak{g}\right)
\]
\end{lemma}

The proof is elementary. Observe that, though $t$ is bounded by the dimension, $k\geq2$
has no restrictions. However, the sum is finite, for the spaces $\mathbf{E}_{c,1}^{t,k_{1},..,k_{r}}$  are zero
for almost any choice $\left( k_{1},.., k_{r}\right)$.\newline Given the Lie algebra  $\frak{g}=\left(\mathbb{C}^{n},\mu_{0}\right)$, we have the associated graduation  
$\frak{gr}\left(  \frak{g}\right)  =\sum_{i=1}^{n\left(  \frak{g}\right)
}\frak{g}_{i},$ where $\frak{g}_{i}=\frac{C^{i-1}\frak{g}}{C^{i}\frak{g}}$ and
$n\left(  \frak{g}\right)  $ is the nilindex of $\frak{g}.$
Independently
of $\frak{g}$ being naturally graded or not, any
vector $X$ has a fixed position in one of the graduation blocks.

\begin{remark}
The study of the central extensions which preserve a graduation is reduced to the study of the position of the adjoined vector $X_{n+1}$. Note that in this sense the cocycles
$\varphi_{ij}\in H_{k}^{2,t}\left(  \frak{g},\mathbb{C}\right)  $ codify 
this information.
\end{remark}

\section{The centralizer property for naturally graded Lie algebras}

For the filiform Lie algebras, Vergne proved that there exist, up to
isomorphism, only two clases of naturally graded Lie algebras, $L_{n}$ and
$Q_{n}$, for which the second only exists in even dimension. Now the model $Q_{n}$ has an interesting structural property that explains its nonexistence in odd dimension : for $p\geq \left[\frac{n-1}{2}\right]$, where [ ] denotes the integer part function, the centralizer $C_{Q_{n}}\left(C^{p}Q_{n}\right)$ in $Q_{n}$ contains the ideal $C^{p}Q_{n}$, while for $q\leq p$ we have $C_{Q_{n}}\left( C^{p}Q_{n}\right) \varsupsetneq C^{p}Q_{n}$. This property gives an estimation of how far the algebra $Q_{n}$ is from having an abelian commutator algebra, as happens for $L_{n}$. This fact suggestts to study the Lie algebras $\frak{g}$ which satisfy the following property :
\begin{eqnarray}
C_{\frak{g}}\left( C^{p}\frak{g}\right)  &\supseteq &C^{p}\frak{g}\text{,}%
\frak{\;}p\geq \left[ \frac{n\left( \frak{g}\right) }{2}\right]   \label{P}
\\
C_{\frak{g}}\left( C^{p}\frak{g}\right)  &\varsupsetneq &C^{p}\frak{g},\;p<%
\left[ \frac{n\left( \frak{g}\right) }{2}\right]   \notag
\end{eqnarray}
where $n\left( \frak{g}\right) $\bigskip\ is the nilpotence class ( or
nilindex ) of $\frak{g}$. We will call $\left(P\right)$ the centralizer property.\bigskip  

Now, a detailed analysis shows that $\left(P\right)$ can be satisfied in two manners :

\begin{enumerate}
\item  There exist $X,Y\in C^{\left[ \frac{n\left( \frak{g}\right) }{2}%
\right] -1}\frak{g-}$ $C^{\left[ \frac{n\left( \frak{g}\right) }{2}\right] }%
\frak{g}$ such that $\left[ X,Y\right] \neq 0,$

\item  if there are $X,Y\in C^{\left[ \frac{n\left( \frak{g}\right) }{2}%
\right] -1}\frak{g}$ with $\left[ X,Y\right] \neq 0$, then either $X\in C^{%
\left[ \frac{n\left( \frak{g}\right) }{2}\right] }\frak{g}$ or $Y\in C^{%
\left[ \frac{n\left( \frak{g}\right) }{2}\right] }\frak{g}$.
\end{enumerate}

We have also to distinguish between two classes of algebras :

\begin{definition}
A naturally graded Lie algebra satisfying the centralizer property through condition $1.$ ( thus not verifying 2.) is called a $\left(P1\right)$-algebra.\newline A naturally graded Lie algebra satisfying the centralizer property through condition $2.$ ( thus not verifying 1.) is called a $\left(P2\right)$-algebra.
\end{definition}  

\begin{remark}
It seems that condition 1) is less natural than 2), for it implies a nonzero bracket
in a specific graduation block of $\frak{g}$. However, we will see that even
those algebras satisfying 1) are obtained in a ''natural'' manner.
\end{remark}

\subsection{$\left(P1\right)$-algebras}

Let $\frak{g}$ be a complex semisimple Lie algebra of finite dimension,
$\frak{h}$ a Cartan subalgebra, $\Phi$ a root system associated to $\frak{h}$
and $\Delta$ a basis of simple roots. We introduce a partial ordering in
$\frak{h}^{\ast}$ relative to which the elements are called positive if their
are linear combinations of simple roots with nonnegative coefficients [8].
Thus, respect to this ordering, we have
\[
\Phi=\Phi^{+}\cup\Phi^{-}%
\]
Recall that the subalgebra $\frak{h}$ induces the Cartan decomposition
\[
\frak{g}=\frak{h}+\coprod_{\alpha\in\Phi}L_{\alpha}%
\]
into weight spaces.\newline
Let $ht:\Phi\rightarrow\mathbb{Z}^{+}$ be the height function. We can define the sets%
\[
\Delta\left(  k\right)  =\left\{  \alpha\in\Phi^{+}\;|\;ht\left(  \alpha\right)  =k\right\}
\]
These sets are of importance, as they give us a natural graduation of the
nilradical of a standard Borel subalgebra $\frak{b}\left(\Delta\right)$.

\begin{theorem}
Let $\frak{n}$ be the nilradical of a standard Borel subalgebra $\frak{b}%
\left(  \Delta\right)  $ of a complex simple Lie algebra distinct from $G_{2}%
$. Then $\frak{n}$ satisfies $\left(P\right)$ and $1)$.
\end{theorem}

The proof is an immediate consequence of the following result :

\begin{proposition}
Let $\frak{n}$ be the nilradical of a standard Borel subalgebra $\frak{b}%
\left(  \Delta\right)  $ of a complex simple Lie algebra distinct from $G_{2}%
$. Let $p=ht\left(  \delta\right)  $ be the height of the maximal root. Then
there exist roots $\alpha,\beta$  whose height is $[\frac{ht\left(  \delta\right)  }{2}]$ such that $\alpha+\beta$ is a positive root.
\end{proposition}

\begin{proof}
\begin{enumerate}
\item $\frak{g}=A_{l}$ :

\begin{enumerate}
\item $l=2q$. Let $\delta$ be the maximal root. For $1\leq t\leq q-2$ we have
\[
\left(  \delta-\alpha_{2q}-\alpha_{2q-1}-..-\alpha_{2q-t},\alpha
_{2q-t-1}\right)  =1
\]
which proves that $\omega_{1}=\alpha_{1}+..+\alpha_{q}$ is a root. In the same
way it is seen that $\omega_{2}=\alpha_{q+1}+..+\alpha_{2q}$ is also a root.
We have $ht\left(  \omega_{1}\right)  =ht\left(  \omega_{2}\right)  =q$, thus
$\omega_{1},\omega_{2}\in\Delta\left(  q\right)  $.

\item $l=2q+1$. Reasoning as before, it follows that
\[
\omega_{1}=\alpha_{1}+..+\alpha_{q},\;\omega_{2}=\alpha_{q+1}+..+\alpha
_{2q}\in\Delta\left(  q\right)
\]
and $\omega_{1}+\omega_{2}=\delta-\alpha_{n}\in\Phi$, as $\alpha_{n}$ is a
particular root.
\end{enumerate}

\item $\frak{g}=B_{l}$ : For $1\leq t\leq l-2$ we have
\[
\left(  \delta-\alpha_{2}-..-\alpha_{t},\alpha_{t+1}\right)  >0
\]
so that $\omega_{1}^{\prime}=\delta-\left(  \alpha_{2}+..+\alpha_{l-1}\right)
$. Now $\left(  \omega_{1}^{\prime},\alpha_{1}\right)  >0$ and $\left(
\omega_{1}^{\prime}-\alpha_{1},\alpha_{2}\right)  >0$, thus
\[
\omega_{1}=\alpha_{3}+..+\alpha_{l-1}+2\alpha_{l}\in\Delta\left(  l-1\right)
=\Delta\left(  \left[  \frac{ht\left(  \delta\right)  }{2}\right]  \right)
\]
Considering the $\alpha_{t}$-string through $\left(  \alpha_{1}+..+\alpha
_{t-1}\right)  $ $\left(  2\leq t\leq l-1\right)  $ we obtain that
\[
\omega_{2}=\alpha_{1}+..+\alpha_{l-1}\in\Delta\left(  l-1\right)
\]
and $\omega_{1}+\omega_{2}=\delta-\alpha_{2}\in\Phi$.

\item $\frak{g}=C_{l}$ : Consider the maximal root $\delta=2\alpha
_{1}+..+2\alpha_{l-1}+\alpha_{l}$ and the particular root $\alpha_{1}$. \ Now
\[
\left(  \delta_{1},\alpha_{j}\right)  =\left\{
\begin{array}
[c]{l}%
1\text{ \ for }j=1\Rightarrow\omega_{1}=\delta_{1}-\alpha_{1}\in\Delta\left(
l-1\right)  \\
1\text{ \ for }j=l\;\Rightarrow\omega_{2}=\delta_{1}-\alpha_{l}\in
\Delta\left(  l-1\right)
\end{array}
\right.
\]
and $\omega_{1}+\omega_{2}=\delta-\alpha_{1}\in\Delta\left(  2l-2\right)  $.

\item $\frak{g}=D_{l}$ : As before, we have $\delta_{1}-\alpha_{1},\delta
_{1}-\alpha_{l}\in\Phi$. Considering the $\alpha_{2}$-string through
$\delta_{1}-\alpha_{1}$ and the $\alpha_{l-1}$-sting through $\delta
_{1}-\alpha_{l}$ we obtain $\omega_{1}=\delta_{1}-\alpha_{1}-\alpha_{2}%
,\omega_{2}=\delta_{1}-\alpha_{l}-\alpha_{l-1}\in\Delta\left(  l-2\right)
\,$\ and $\omega_{1}+\omega_{2}=\delta-\alpha_{2}\in\Delta\left(  2l-4\right)
$.

\item $\frak{g}=E_{6}$ : $ht\left(  \delta\right)  =11$%
\begin{align*}
\omega_{1} &  =\alpha_{1}+\alpha_{3}+\alpha_{4}+\alpha_{5}+\alpha_{6}\in
\Delta\left(  5\right)  \\
\omega_{2} &  =\alpha_{2}+\alpha_{3}+2\alpha_{4}+\alpha_{5}\in\Delta\left(
5\right)
\end{align*}
where $\omega_{1}+\omega_{2}=\delta-\alpha_{2}$.

\item $\frak{g}=E_{7}$ : $ht\left(  \delta\right)  =17$%
\begin{align*}
\omega_{1} &  =\alpha_{1}+\alpha_{2}+2\alpha_{3}+2\alpha_{4}+\alpha_{5}%
+\alpha_{6}\in\Delta\left(  8\right)  \\
\omega_{2} &  =\alpha_{2}+\alpha_{3}+2\alpha_{4}+2\alpha_{5}+\alpha_{6}%
+\alpha_{7}\in\Delta\left(  8\right)
\end{align*}
$\omega_{1}+\omega_{2}=\delta-\alpha_{1}.$

\item $\frak{g}=E_{8}$ : $ht\left(  \delta\right)  =29$%
\begin{align*}
\omega_{1} &  =\alpha_{1}+2\alpha_{2}+2\alpha_{3}+3\alpha_{4}+3\alpha
_{5}+2\alpha_{6}+\alpha_{7}\in\Delta\left(  14\right)  \\
\omega_{2} &  =\alpha_{1}+\alpha_{2}+2\alpha_{3}+3\alpha_{4}+2\alpha
_{5}+2\alpha_{6}+2\alpha_{7}+\alpha_{8}\in\Delta\left(  14\right)
\end{align*}
$\omega_{1}+\omega_{2}=\delta-\alpha_{8}$.

\item $\frak{g}=F_{4}$ : $ht\left(  \delta\right)  =11$%
\begin{align*}
\omega_{1} &  =\alpha_{1}+2\alpha_{2}+2\alpha_{3}\in\Delta\left(  5\right)  \\
\omega_{2} &  =\alpha_{2}+2\alpha_{3}+2\alpha_{4}\in\Delta\left(  5\right)
\end{align*}
$\omega_{1}+\omega_{2}=\delta-\alpha_{1}$.
\end{enumerate}
\end{proof}

\begin{remark}
We see that the classical theory provides a lot of naturally graded Lie algebras satisfying the centralizer property. However, it is usually unconvenient to manipulate these algebras, because of the great difference between its dimension and nilpotence class : the first is too high in comparison with the last.
\end{remark}

\begin{note}
Unless otherwise stated, whenever we speak in future about Lie algebras $\frak{g}$ satisfyng $\left(P\right)$, we will understand that $\frak{g}$ is naturally graded and satisfies condition $2)$ and not $1)$.
\end{note}

\subsection{$\left(P2\right)$-algebras}

With the conventions adopted, it follows immediately from the preceding remarks :

\begin{proposition}
A $\left(P2\right)$ Lie algebra $\frak{g}$ is filiform if and only if $\frak{g}\simeq Q_{n}$ for $n\geq 6$.
\end{proposition}

Recall that the set of characteristic sequences for a nilpotent Lie algebra $\frak{g}$ is ordered lexicographically. As $\left(n-1,1\right)$ is the maximum of this set, it is natural to begin our study with its immediate successor, i.e, the sequence $\left(n-2,1,1\right)$, where $n=dim\left(\frak{g}\right)$. If $X_{1}$ is a characteristic vector, then we can find a basis $\left\{  X_{1},...,X_{n}\right\}  $ such that $\left[  X_{1},X_{i}\right]  =X_{i+1},\;2\leq
i\leq n-2.$ Let $\omega_{1}$ be the vector of the dual base which corresponds
to $X_{1}$. It follows from the brackets that the exterior product $\omega
_{1}\wedge\omega_{i-1}$ is a summand of the differential form $d\omega_{i}$. If
$\frak{g}$ satisfies $\left(P2\right)$, then the existence of a vector $X_{j}$
$\left(  j\geq2\right)  $ with $\left[  X_{k+1},X_{j}\right]  \neq0$ and
$\left[  X_{j},X_{j+t}\right]  =0,\;$for all $t\geq1$ suffices.

\begin{proposition}
Let $\frak{g}$ be a $n$-dimensional Lie algebra with characteristic sequence
$\left(  n-2,1,1\right)  $. If $X_{n}\in\frak{g}_{2t}$, $1\leq t\leq\frac
{n-2}{2}$ then $\frak{g}$ is not naturally graded.
\end{proposition}

\begin{proof}
We have, for any $t$, $\frak{g}_{2t}=\frac{C^{2t-1}\frak{g}}{C^{2t}\frak{g}}$.
If $X_{n}\in\frak{g}_{2t}$ we have the brackets
\begin{align*}
\left[  X_{2},X_{2t}\right]   &  =\lambda_{2t-2}X_{2t+1}+\mu_{1}X_{n}\\
\left[  X_{j},X_{2t-j+2}\right]   &  =\lambda_{2t-2}^{j}X_{2t+1}+\mu_{1}%
^{t}X_{n}%
\end{align*}
Applying the adjoint operator $ad\left(  X_{1}\right)  $ we obtain the
condition :%
\[
\lbrack X_{1},[X_{t},X_{t+1}]]+\left[  X_{t},X_{t+2}\right]  =0
\]
Now $ad\left(  X_{t}\right)  \left(  \left\langle X_{1},..,X_{t+1}%
\right\rangle _{\mathbb{C}}\right)  \cap\left\langle X_{n}\right\rangle
=\{0\}$, so that $\mu_{1}^{t}=0$ for all $t$. On the other hand, the previous
condition implies $\mu_{1}=\left(  -1\right)  ^{t-1}\mu_{1}^{t}\;\forall t$,
so $X_{n}\notin C^{1}\frak{g}$, contradiction with the assumption.
\end{proof}

\begin{remark}
As a consequence of the previous result, the position of a vector $X_{n}$ is
never optimal in an even indexed graduation block. For this reason, it is convenient to introduce the following convention : for a $n$-dimensional nilpotent Lie algebra $\frak{g}$ with basis $\left\{X_{1},..,X_{n}\right\}$ we say that the vector $X_{n}$ has depth $k$, noted $h\left(X_{n}\right)=k\; \left(1\leq k\leq \left[\frac{n-3}{2}\right]\right)$, if $X_{n}\in \frak{g}_{2k+1}$.\newline In fact, this definition can be extended to any vector of $\frak{g}$. Observe that there are fractional depths.
\end{remark}

\begin{theorem}
There do not exist even dimensional $\left(P2\right)$ Lie algebras $\frak{g}$ of characteristic sequence $\left(  2m-2,1,1\right)$.
\end{theorem}

\begin{proof}
The characteristic sequence imposes the existence of a basis $\left\{
X_{1},..,X_{n}\right\}  $ such that
\[
adX_{1}\left(  X_{i}\right)  =X_{i+1},\;2\leq i\leq2m-2
\]
thus we have
\[
C^{k}\frak{g}\supset\left\langle X_{k+2},..,X_{2m-1}\right\rangle ,\;1\leq
k\leq2m-3
\]
The central descending sequence induces the following relations for the
associated graduation
\[
\frak{g}_{1}\supset\left\langle X_{1},X_{2}\right\rangle ,\;\frak{g}%
_{k}\supset\left\langle X_{k+1}\right\rangle ,\;2\leq k\leq2m-2
\]
If $\left(P2\right)$ is satisfied, there exist two nonzero
vectors $X,Y\in C^{m-2}\frak{g}$ such that $\left[  X,Y\right]  \neq0$.
Without loss of generality we can suppose $X=X_{m},\;Y=X_{m+t}$ for $t\geq1$.
Then
\[
\left[  X_{m},X_{m+t}\right]  \in\left[  \frak{g}_{m-1},\frak{g}%
_{m+t-1}\right]  \subset\frak{g}_{2m-2+t}=\{0\},\;t\geq1
\]
This shows that the unique admissible graduation block for  $X_{2m}$ is $\frak{g}_{m-1}$. Suppose therefore that $X_{2m}\in \frak{g}_{m-1}$. Then $\left(P2\right)$ implies $\left[ X_{2m}, X_{m}\right]=\lambda X_{2m-1}$ for a nonzero value $\lambda$; moreover, $m$ must be  even, $m=2r$. As $X_{2m}$ belongs to the commutator algebra, there exist two indexes $i,j\geq 1$ with $i+j=m-1$ and a pair of vectors $X_{i+1}\in \frak{g}_{i}$, $X_{j+1}\in \frak{g}_{j}$ such that $\left[ X_{i+1}, X_{j+1}\right]=\lambda_{i,j} X_{2m}$, where $\lambda_{i,j}$ is nonzero. Let $\left( i_{0}, j_{0}\right)$ be the minimal pair with this property; it is not difficult to see that it is $\left( 2,2r-3\right)$. Then the associated differential form to the vector $X_{2m}$ is of the following type :
\[
d\omega_{2m}= \sum_{t\geq 0} \lambda_{1+t,2r-4-t} \omega_{2+t}\wedge\omega_{2r-3-t}
\]
On the other hand 
\[
d\omega_{2m-1}= \sum_{s\geq 0} \alpha_{t} \omega_{2+s}\wedge\omega_{2m-3-t}+\lambda \omega_{m}\wedge\omega{2m}
\]
It is immediate to verify the nonexistence of nonzero coefficients $\lambda_{1+t,2r-4-t}$ such that the previous forms satisfy simultaneously
\[
d^{2}\omega_{2m-1}=d^{2}\omega_{2m}=0
\]
\end{proof}

\begin{remark}
The obstruction for the even dimensional is the same as the one observed in
the analysis of filiform algebras. In this sense, the ( odd ) dimensional Lie
algebras which verify $\left(P2\right)$ will play the same role
as $Q_{n}$ does for the filiform algebras.
\end{remark}

Now we approach the classification problem : to obtain all nilpotent Lie algebras that satisfy $\left(P2\right)$ and whose characteristic sequence is $\left( 2m-1,1,1\right)$. To avoid trivial cases, the algebras are supposed to be nonsplit.

\bigskip For $m\geq4$ let $\frak{s}_{m}$ be $2m$-dimensional Lie algebra whose structural equations over the
basis $\left\{  \omega_{1},..,\omega_{2m-1},\omega_{2m+1}\right\}  $ of
$\left(  \mathbb{C}^{2m}\right)  ^{\ast}$ are
\begin{align*}
d\omega_{1}  &  =d\omega_{2}=0\\
d\omega_{j}  &  =\omega_{1}\wedge\omega_{j-1},\;3\leq j\leq2m-3\\
d\omega_{2m-2}  &  =\omega_{1}\wedge\omega_{2m-3}+\sum_{j=2}^{\left[
\frac{2m+1}{2}\right]  -1}\left(  -1\right)  ^{j}\omega_{j}\wedge
\omega_{2m-1-j}\\
d\omega_{2m-1}  &  =\omega_{1}\wedge\omega_{2m-2}-\left(  m-2\right)
\omega_{2}\wedge\omega_{2m+1}+\sum_{j=2}^{m-1}\left(  -1\right)  ^{j}\left(
m-i\right)  \omega_{j}\wedge\omega_{2m-j}\\
d\omega_{2m+1}  &  =\sum_{j=2}^{\left[  \frac{2m+1}{2}\right]  -1}\left(
-1\right)  ^{j}\omega_{j}\wedge\omega_{2m-1-j}%
\end{align*}
It follows immediately that $\frak{s}_{m}$ is naturally graded of
characteristic sequence $\left(  2m-2,1,1\right)  $ for any $m\geq4.$

\begin{notation}
Consider $\mathbf{E}_{c,1}\left(\frak{g}\right)$ and let $p\in\left\{n\left(\frak{g}\right), n\left(\frak{g}\right)+1\right\}$. Denote by  $\mathbf{E}_{c,1}\left(\frak{g},p\right)$ the extensions whose nilindex is $p$.
\end{notation}

\begin{proposition}
Let $m\geq 4$ and $e\left(\frak{s}_{m}\right)\in\mathbf{E}_{c,1}^{2m-1}\left(\frak{s}_{m}\right)$. Then $e\left(\frak{s}_{m}\right)$ satisfies $\left(P2\right)$ if and only if it is isomorphic to the Lie algebra $\frak{g}_{\left(  m,m-2\right)  }^{3}$ given by :

\begin{align*}
d\omega_{1}  & =d\omega_{2}=0\\
d\omega_{j}  & =\omega_{1}\wedge\omega_{j-1},\;3\leq j\leq2m-3\\
d\omega_{2m-2}  & =\omega_{1}\wedge\omega_{2m-3}+\sum_{j=2}^{\left[
\frac{2m+1}{2}\right]  -1}\left(  -1\right)  ^{j}\;\omega_{j}\wedge
\omega_{2m-1-j}\\
d\omega_{2m-1}  & =\omega_{1}\wedge\omega_{2m-2}-\left(  m-2\right)
\,\omega_{2}\wedge\omega_{2m+1}+\sum_{j=2}^{m-1}\left(  -1\right)  ^{j}\left(
m-j\right)  \,\ \omega_{j}\wedge\omega_{2m-j}\\
d\omega_{2m}  & =\omega_{1}\wedge\omega_{2m-1}-\left(  m-2\right)
\,\omega_{3}\wedge\omega_{2m+1}+\sum_{j=3}^{m}\frac{\left(  -1\right)
^{j}\left(  j-2\right)  \left(  2m-1-j\right)  }{2}\;\omega_{j}\wedge
\omega_{2m+1-j}\\
d\omega_{2m+1}  & =\sum_{j=2}^{\left[  \frac{2m+1}{2}\right]  -1}\left(
-1\right)  ^{j}\;\omega_{j}\wedge\omega_{2m-1-j}%
\end{align*}
for $m\geq 5$.\newline

If $m=4$ there is an additional extension $\frak{g}_{\left(  4,2\right)  }^{1}$ :

\begin{align*}
d\omega_{1} &  =d\omega_{2}=0\\
d\omega_{j} &  =\omega_{1}\wedge\omega_{j-1},\;3\leq j\leq5\\
d\omega_{6} &  =\omega_{1}\wedge\omega_{5}+\omega_{2}\wedge\omega_{5}%
-\omega_{3}\wedge\omega_{4}\\
d\omega_{7} &  =\omega_{1}\wedge\omega_{6}+2\omega_{2}\wedge\omega_{6}%
-\omega_{3}\wedge\omega_{5}-2\omega_{2}\wedge\omega_{9}\\
d\omega_{8} &  =\omega_{1}\wedge\omega_{7}+\omega_{2}\wedge\omega_{7}%
-\omega_{3}\wedge\omega_{6}+2\omega_{4}\wedge\omega_{5}-2\omega_{3}%
\wedge\omega_{9}\\
d\omega_{9} &  =\omega_{2}\wedge\omega_{5}-\omega_{3}\wedge\omega_{4}%
\end{align*}
\end{proposition}

\begin{proof}
Let $\left\{\omega_{1},..,\omega_{2m},\omega_{2m+1}\right\}$ be a basis of  $e\left(\frak{s}_{m}\right)$ over $\left(\mathbb{C}^{2m+1}\right)^{*}$ and $\left\{X_{1},..,X_{2m+1}\right\}$ its dual basis.
\newline Any central extension is specified by the adjunction of a differential
form $d\omega_{2m}.$ The graduation forces the depth of $X_{2m}$ \ to be
$h\left(  X_{2m}\right)  =m-2.$ So $d\omega_{2m}$ is of the following type
\[
d\omega_{2m}=\omega_{1}\wedge\omega_{2m-1}+\sum_{j=2}^{m-1}\varphi
_{j,2m-j}\,\omega_{j}\wedge\omega_{2m+1-j}+\varphi_{2,2m+1}\omega_{2}%
\wedge\omega_{2m+1}%
\]
where $\varphi_{j,2m-j}\in H_{2}^{2,m-1}\left(  \frak{s}_{m}\right)  $ \ for
$j=2,..,\left[  \frac{2m+1}{2}\right]  $ \ and \ $\varphi_{2,2m+1}\in
H_{4}^{2,m-1}\left( \frak{s}_{m},\mathbb{C}\right)  $. The structure of $\frak{s}_{m}$
implies $0\neq\varphi_{2,2m+1}.$ Moreover, the following relations hold
\[
\left(  j-2\right)  \left(  2m-1-j\right)  \varphi_{3,2m-2}+2\left(
-1\right)  ^{j}\left(  m-2\right)  \varphi_{j,2m-j}=0,\;j=4,..,\left[
\frac{2m+1}{2}\right]
\]%
\[
\frac{m}{2}\left(  m-3\right)  \varphi_{m-1,m+1}+\varphi_{m-1,m+2}=0
\]
from which we deduce, by the structure of $\frak{s}_{m},$ that $\varphi_{2,2m-1}=0$.
Observe in particular that the nullity of this cocycle implies the existence
of a unique extension. Through an elementary change of basis it follows that
this extension is isomorphic to $\frak{g}_{\left(  m,m-1\right)  }^{3}$ \ for
$m\geq5$.\newline
An algebra $e_{c,1}\left(\frak{s}_{4}\right)  \in\mathbf{E}_{c,1}\left(
\frak{s}_{4}\right)$ \ \ is determined by the adjunction of a
differential form $d\omega_{8}$. As the nilindex $p=8$ is fixed, this implies
that $h\left(  X_{8}\right)  =3.$ Then this form must be of the following type
:%
\[
d\omega_{8}=\omega_{1}\wedge\omega_{7}+\varphi_{27}\omega_{2}\wedge\omega
_{7}+\varphi_{36}\omega_{3}\wedge\omega_{6}+\varphi_{45}\omega_{4}\wedge
\omega_{5}+\varphi_{29}\omega_{2}\wedge\omega_{9}%
\]
where $\varphi_{27},\varphi_{36},\varphi_{45}\in H_{2}^{2,3}\left(
\frak{s}_{4},\mathbb{C}\right)  $ and $0\neq\varphi_{29}\in H_{4}^{2,3}\left(
\frak{s}_{4},\mathbb{C}\right)  .$ The determinant cocycle is $\varphi_{27}$ : if it
is nonzero we obtain
\[
X_{2}\wedge X_{7}-X_{3}\wedge X_{6},2X_{2}\wedge X_{7}+X_{4}\wedge X_{5}%
\in\Omega
\]
and otherwise
\[
3X_{3}\wedge X_{6}+2X_{4}\wedge X_{5}\in\Omega
\]
Thus there are two nonequivalent extensions, the first being isomorphic to
$\frak{g}_{\left(  4,2\right)  }^{1}$ and the second to $\frak{g}_{\left(
4,2\right)  }^{2}$.
\end{proof}

\begin{theorem}
If $h\left(  X_{2m+1}\right)  =t,\;t\neq m-2,$ any naturally graded Lie algebra $\frak{g}$ \ with characteristic
sequence $\left(  2m-1,1,1\right)$ that satisfies $\left(P2\right)$ is a central extension of either $Q_{2m-1}$ or $L_{2m-1}$.
\end{theorem}

\begin{proof}
Let $h\left(  X_{2m+1}\right)  =t,\;t\neq m-2.$ We define the cocycle
$\varphi_{j2m+1}\in H^{2}\left(  \frak{g},\frak{g}\right)  $ by
\[
\varphi_{j2m+1}\left(  X_{j},X_{2m+1}\right)  =\left\{
\begin{array}
[c]{r}%
\alpha_{j2m+1}X_{j+1+2t}\;\;\text{if \ \ }j+2t\leq2m-1\\
0\text{\ \ if\ }j+2t>2m-1
\end{array}
\right.
\]
The action of the adjoint operator $ad\left(  X_{1}\right)  $ \ implies the
conditions
\[
\alpha_{22m+1}=-\alpha_{32m+1}=....=\left(  -1\right)  ^{k_{0}}\alpha
_{k_{0}2m+1}%
\]
where $k_{0}$ is the last value for which $\varphi_{j2m+1}$ is nonzero.
Moreover,
\[
d\omega_{2m}=\omega_{1}\wedge\omega_{2m-1}+\sum_{i,j}\varphi_{ij}\,\omega
_{i}\wedge\omega_{j}+\alpha_{2,2m+1}\omega_{2}\wedge\omega_{2m+1}%
\]
where $\varphi_{ij}\in H_{2}^{2,m-1}\left(  \frak{g},\mathbb{C}\right)  .$ The
Jacobi condition $d\left(  d\omega_{2m}\right)  =0$ implies $\alpha
_{2,2m+1}=0,$ so that the cocycle $\varphi_{j2m+1}$ is identically zero. Thus
the vector $X_{2m+1}$ is central and the factor algebra $\frac{\frak{g}%
}{\left\langle X_{2m+1}\right\rangle }$ \ is naturally graded and filiform,
isomorphic to $Q_{2m-1}$ if $t\neq m-1$ and isomorphic to $L_{2m-1}$ if $t=m-1.$
\end{proof}

\begin{remark}
We commented the existence of vectors having fractionary depth, according to the definition given before. To cover all cases, it must be shown that for these depths there do not exist extensions which satisfy the required conditions.
\end{remark}

\begin{proposition}
For $m\geq4$ and $q\equiv0\left(  \operatorname{mod}\;2\right)  $
\[
\mathbf{E}_{c,1}^{\frac{q+1}{2},2}\left(  Q_{2m-1},2m-1\right)=0\text{ }%
\]
\end{proposition}

\begin{proof}
An extension $e\left(  Q_{2m-1}\right)  \in$ $\mathbf{E}_{c,1}^{\frac{q+1}%
{2},2}\left(  Q_{2m-1},p\right)$ is determined by the
cocycles of \ the space $H_{2}^{2,\frac{q+1}{2}}\left(  Q_{2m-1}%
,\mathbb{C}\right)  .$ Then the differential form $d\omega_{2m+1}$ is of type
\[
d\omega_{2m+1}=\sum_{i,j}\varphi_{ij}\,\omega_{i}\wedge\omega_{j}%
\]
where the indexes $i,j$ satisfy
\[
i+j=q+4
\]
As $q$ is even, let $q=2t.$ The the form $d\omega_{2m+1}$ can be rewritten as
\[
d\omega_{2m+1}=\sum_{j=2}^{t+1}\varphi_{j,4+2t-j}\,\omega_{j}\wedge
\omega_{4+2t-j}%
\]
It is trivial to verify that the equations
\[
\varphi_{2,2+2t}+\left(  -1\right)  ^{j}\varphi_{j,4+2t-j}=0,\ 3\leq j\leq
t+1
\]
are satisfied. This allows us to take a common factor, so that
\[
d\omega_{2m+1}=\varphi_{2,2+2t}d\varpi_{t}%
\]
where this form is easily proven to be nonclosed. So we deduce the nonexistence of naturally graded with the required
nilindex in $\mathbf{E}_{c,1}^{\frac{q+1}{2},2}\left(  Q_{2m-1}\right)  .$
\end{proof}

\begin{proposition}
For $m\geq4$ and $q\equiv0\left(  \operatorname{mod}\;2\right)  $
\[
\mathbf{E}_{c,1}^{\frac{q+1}{2},2}\left(  L_{2m-1},2m-1\right)=0\text{ }%
\]
\end{proposition}

The proof is analogous to the previous case.

\begin{corollary}
For $m\geq5, r\geq 2$ and $k_{1},..,k_{r}\in\mathbb{Z}^{+}$%
\[
\mathbf{E}_{c,1}^{\frac{q+1}{2},k_{1},..,k_{r}}\left(  Q_{2m-1},2m-1\right)
=\mathbf{E}_{c,1}^{\frac{q+1}{2},k_{1},..,k_{r}%
}\left(  L_{2m-1},2m-1\right)=\{0\}
\]
\end{corollary}

\begin{theorem}[Classification of $\left(P2\right)$-algebras with ch.s.
$\left(  2m-1,1,1\right)$]
A naturally graded Lie algebra $\frak{g}$ with characteristic sequence $\left(2m-1,1,1\right)$ is a $\left(P2\right)$-algebra
if and only if it is isomorphic to one of the following models :

\begin{enumerate}
\item $\frak{g}_{\left(  4,2\right)  }^{1}:$%
\begin{align*}
d\omega_{1}  &  =d\omega_{2}=0\\
d\omega_{j}  &  =\omega_{1}\wedge\omega_{j},\;j=3,4,5\\
d\omega_{6}  &  =\omega_{1}\wedge\omega_{5}+\omega_{2}\wedge\omega_{5}%
-\omega_{2}\wedge\omega_{4}\\
d\omega_{7}  &  =\omega_{1}\wedge\omega_{6}+2\omega_{2}\wedge\omega_{6}%
-\omega_{3}\wedge\omega_{5}-2\omega_{2}\wedge\omega_{9}\\
d\omega_{8}  &  =\omega_{1}\wedge\omega_{7}+\omega_{2}\wedge\omega_{7}%
+\omega_{3}\wedge\omega_{6}-2\omega_{4}\wedge\omega_{5}-2\omega_{3}%
\wedge\omega_{9}\\
d\omega_{9}  &  =\omega_{2}\wedge\omega_{5}-\omega_{3}\wedge\omega_{4}%
\end{align*}

\item $\frak{g}_{\left(  m,t\right)  }^{2}\;\left(  1\leq t\leq m-2\right)  :$%
\begin{align*}
d\omega_{1}  &  =d\omega_{2}=0\\
d\omega_{j}  &  =\omega_{1}\wedge\omega_{j-1},\;3\leq j\leq2m-1\\
d\omega_{2m}  &  =\omega_{1}\wedge\omega_{2m-1}+\sum_{j=2}^{\left[
\frac{2m+1}{2}\right]  }\left(  -1\right)  ^{j}\;\omega_{j}\wedge
\omega_{2m+1-j}\\
d\omega_{2m+1}  &  =\sum_{j=2}^{t+1}\left(  -1\right)  ^{j}\;\omega_{j}%
\wedge\omega_{3-j+2t}%
\end{align*}

\item $\frak{g}_{\left(  m,m-2\right)  }^{3}:$%
\begin{align*}
d\omega_{1}  &  =d\omega_{2}=0\\
d\omega_{j}  &  =\omega_{1}\wedge\omega_{j-1},\;3\leq j\leq2m-3\\
d\omega_{2m-2}  &  =\omega_{1}\wedge\omega_{2m-3}+\sum_{j=2}^{\left[
\frac{2m+1}{2}\right]  -1}\left(  -1\right)  ^{j}\;\omega_{j}\wedge
\omega_{2m-1-j}\\
d\omega_{2m-1}  &  =\omega_{1}\wedge\omega_{2m-2}+\left(  m-2\right)
\,\omega_{2}\wedge\omega_{2m+1}+\sum_{j=2}^{m-1}\left(  -1\right)  ^{j}\left(
m-j\right)  \,\ \omega_{j}\wedge\omega_{2m-j}\\
d\omega_{2m}  &  =\omega_{1}\wedge\omega_{2m-1}+\left(  m-2\right)
\,\omega_{3}\wedge\omega_{2m+1}+\sum_{j=3}^{m}\frac{\left(  -1\right)
^{j}\left(  j-2\right)  \left(  2m-1-j\right)  }{2}\;\omega_{j}\wedge
\omega_{2m+1-j}\\
d\omega_{2m+1}  &  =\sum_{j=2}^{\left[  \frac{2m+1}{2}\right]  -1}\left(
-1\right)  ^{j}\;\omega_{j}\wedge\omega_{2m-1-j}%
\end{align*}

\item $\frak{g}_{\left(  m,m-1\right)  }^{4}:$%
\begin{align*}
d\omega_{1}  &  =d\omega_{2}=0\\
d\omega_{j}  &  =\omega_{1}\wedge\omega_{j-.1},\;3\leq j\leq2m\\
d\omega_{2m+1}  &  =\sum_{j=2}^{\left[  \frac{2m+1}{2}\right]  }\left(
-1\right)  ^{j}\omega_{j}\wedge\omega_{2m+1-j}\;
\end{align*}
Moreover, these algebras are pairwise non isomorphic.
\end{enumerate}
\end{theorem}

The proof will be a consequence of the next results :

\begin{lemma}
For $m\geq 4$ the following equations hold

\begin{enumerate}
\item $\mathbf{E}_{c,1}\left(  Q_{2m-1},2m-1\right)=\sum_{t=1}^{t-2}\mathbf{E}_{c,1}^{t,2}\left(  Q_{2m-1},2m-1\right)$

\item $\mathbf{E}_{c,1}\left(  L_{2m-1},2m-1\right)=\mathbf{E}_{c,1}^{m-1,2}\left(  L_{2m-1},2m-1\right)$.
\end{enumerate}
\end{lemma}

\begin{proof}
It is not difficult to see that if $k\neq2$, then
\[
\mathbf{E}_{c,1}^{t,k}\left(  \frak{g},2m-1\right)=0,\;\frak{g}=L_{n}\text{ or \ }Q_{n}%
\]
For $k=2$ and any of the nongiven $t$'s \ the nonexistence of naturally graded extensions with
the required nilindex is routine. The remaining cases are a direct consequence
of the previous results.
\end{proof}

\begin{proposition}
For $m\geq 4$ \ any extension $\frak{g}'\in\mathbf{E}_{c,1}^{t,2}\left(  Q_{2m-1},2m-1\right)$ is isomorphic to $\frak{g}_{m,t}^{2}$ if\; $1\leq t\leq m-2$. For $t<1$ and $t\geq m-1$\;  $\mathbf{E}_{c,1}^{t,2}\left(  Q_{2m-1},2m-1\right)=0$.
\end{proposition}

\begin{proof}
For $1\leq t\leq m-2$ \ the cocycles $\varphi_{ij}\in H_{2}%
^{2,t}\left(  Q_{2m-1},\mathbb{C}\right)  $ must satisfy the relation
$\ i+j=2t+3.$ It is immediate to verify that this space is generated by the
cocycles
\[
\varphi_{2,2t+1},\varphi_{3,2t},...,\varphi_{t+1,t+2}%
\]
subjected to the relations
\[
\varphi_{2,2t+1}+\left(  -1\right)  ^{j-1}\varphi_{j,2t+3-j}=0,\;j=3,..,t+1
\]
If $\left\{  X_{1},..,X_{2m}\right\}  $ is the dual base of $\left\{
\omega_{1},..,\omega_{2m}\right\}  $, we have
\[
X_{2,2t+1},X_{3,2t},...,X_{t+1,t+2}\in Ker\,\lambda
\]
and
\[
X_{2,2t+1}+\left(  -1\right)  ^{j-1}X_{j,2t+3-j}\in\Omega,\;j=3,..,t+1
\]
$\;\;\;$\newline Thus there is, for any $t,$ only one extension, which is
isomorphic to $\frak{g}_{\left(  m,t\right)  }^{2}.$ For the remaining values
of $t$ it is easy to see that $Q_{2m-1}$ does not admit naturally graded
extensions with the prescribed characteritic sequence.
\end{proof}

\begin{proposition}
For $m\geq 4$ \ any extension $\frak{g}'\in\mathbf{E}_{c,1}^{t,2}\left(  L_{2m-1},2m-1\right)$ is isomorphic to $\frak{g}_{m,m-1}^{4}$ if $t=m-1$. For $t\neq m-1$  $\mathbf{E}_{c,1}^{t,2}\left(  L_{2m-1},2m-1\right)=0$.
\end{proposition}

\begin{proof}
Similarly to the previous case we have
\[
X_{j}\wedge X_{2m-j}\in Ker\,\lambda,\;j=2,..,\left[  \frac{2m+1}{2}\right]
\]
and
\[
X_{2}\wedge X_{2m-2}+\left(  -1\right)  ^{j-1}X_{j}\wedge X_{2m-j}\in
\Omega,\;j=3,..,\left[  \frac{2m+1}{2}\right]
\]
so that there exists a unique extension, isomorphic to $\frak{g}_{\left(
m,m-1\right)  }^{4}.$
\end{proof}

\begin{proof}[\textbf{Proof of theorem 4}]
$\Longrightarrow)$ We can suppose $m\geq5,$
as we have studied the case $m=4$ \ before. We know that if the depth of the
vector $X_{2m+1}$ is $h(X_{2m+1})=t\in\mathbb{Z},\;t\neq m-2$ the factor
algebra $\frac{\frak{g}}{\left\langle X_{2m+1}\right\rangle }$ is naturally
graded and filiform, thus $\frak{g}$ \ is a central extension of either
$Q_{2m-1}$ or $L_{2m-1}.$ Let also be $\ h\left(  X_{2m+1}\right)  =m-2.$ If
$X_{2m+1}$ is central, we obtain again a central extension of $Q_{2m-1}.$ If
not, then the differential form $d\omega_{2m}$ has a nonzero coefficient
associated to the summand $\omega_{2}\wedge\omega_{2m+1}.$ In this case, the
central element to be taken is $X_{2m},$ and it is not difficult to see that
$\frac{\frak{g}}{\left\langle X_{2m}\right\rangle }$ \ is a naturally graded
Lie algebra isomorphic to $\frak{s}_{m}.$ As the central graded \ extensions of this
algebra which increment the nilindex in one unity are unique, this algebra must
be isomorphic to $\frak{g}_{\left(  m,m-2\right)  }^{3}$.

Finally, for the fractionary depths we have seen the nonexistence of
extensions of this type.\newline $\Longleftarrow)$ It is a trivial
verification that the models satisfy the requirements.
\end{proof}

\section{Classification of $\left(P2\right)$-algebras with characteristic sequence $\left(  2m-1,2,1\right)$}

In this section we use the preceding results to establish a classification of $\left(P2\right)$-algebras when the second entry of the characteristic sequence is increased by one. We will see that, with one exception, these algebras are obtained by considering central extensions of the preceding models.\newline  Let $\left(  2m-1,2,1\right)  $ be the characteristic
sequence of $\frak{g}$ and $X_{1}$ a charateristic vector. Then we can find a
basis $\left\{  X_{1},..,X_{2m+2}\right\}  $ dual to the base $\left\{
\omega_{1},..,\omega_{2m+2}\right\}  $ and such that $\left[  X_{1}%
,X_{i}\right]  =X_{i+1},\;2\leq i\leq2m-1$ and $\left[  X_{1},X_{2m+1}\right]
=X_{2m+2}$.

\begin{remark}
In contrast to the preceding case, we will see that now there are split algebras which admit nonsplit naturally graded central extensions of degree one which satisfy $\left(P2\right)$. This will be justified by the existence of a two dimensional Jordan block for the adjoint operator for a characteristic vector.
\end{remark}

For $m\geq4$, let $\frak{g}_{\left(
m,0\right)  }^{1+k}\;\left(  k=0,1\right)$ be the algebras whose Cartan-Maurer equations over the basis $\left\{\omega_{1},..,\omega_{2m+2}\right\}$ are :
\begin{align*}
d\omega_{1}  &  =d\omega_{2}=0\\
d\omega_{j}  &  =\omega_{1}\wedge\omega_{j-1},\;3\leq j\leq2m-1\\
d\omega_{2m}  &  =\omega_{1}\wedge\omega_{2m-1}+\sum_{j=2}^{\left[
\frac{2m+1}{2}\right]  }\left(  -1\right)  ^{j}\,\omega_{j}\wedge
\omega_{2m+1-j}\\
d\omega_{2m+1}  &  =0\\
d\omega_{2m+2}  &  =\omega_{1}\wedge\omega_{2m+1}+k\omega_{2}\wedge
\omega_{2m+1}%
\end{align*}
which are clearly $\left(  m-1\right)  $-abelian.

\begin{lemma}
For $m\geq 4$, a nosnplit naturally extension of $Q_{n}\oplus\mathbb{C}$ satisfies $\left(P2\right)$ if and only if it is isomorphic to either $\frak{g}_{\left(
m,0\right)  }^{1+k}$ for $k=0,1$.
\end{lemma}

\begin{proof}
Bot the graduation and $\left(P2\right)$ imply that the only cocycles that must be considered are those belonging to the space $H_{2m-1}^{2,1}\left(Q_{2m-1}\oplus\mathbb{C}\right)$. Thus the only cohomology classes that give central extensions with the prescribed conditions are $\varphi_{1,2m+1}$ and $\varphi_{2,2m+1}$. The differential form $d\omega_{2m+2}$ associated to the koined vector $X_{2m+2}$ has the form :
\[
d\omega_{2m+2} =\alpha\varphi_{1,2m+1}\omega_{1}\wedge\omega_{2m+1}+ \beta\varphi_{2,2m+1}\omega_{2}\wedge\omega_{2m+1}
\]
where $\alpha,\beta\in\mathbb{C}$. Clearly $\alpha$ must be nonzero, and by a change of basis we can suppose $\alpha=1$. If $\beta=0$ we obtain $\frak{g}_{\left(
m,0\right)  }^{1}$, while for nonzero $\beta$ we obtain $\frak{g}_{\left(
m,0\right)  }^{2}$.
\end{proof}

\begin{remark}
Observe that with the definition of depth introduced earlier, the vector $X_{2m+2}$ of an algebra $\frak{g}$ of characteristic sequence $\left(2m-1,2,1\right)$ has fractionary depth. This follows directly from it, as the position of this last vector is determined by the one of $X_{2m+1}$, as these vectors form the two dimensional Jordan box for $ad\left(X_{1}\right)$.
\end{remark}

\begin{theorem}
Let $m\geq 4$. If $h\left(X_{2m+2}\right)=\frac{2t+1}{2}$ \;$\left(3\geq t\geq m-2\right)$,\;then a $\left(P2\right)$-algebra of characteristic sequence $\left(2m-1,2,1\right)$ is an extension of $\frak{g}_{\left(m,t\right)}^{2}$ or $\frak{g}_{\left(m,m-1\right)}^{3}$.
\end{theorem}
 
\bigskip The proof is of the same kind as theorem 3.\newline For the lowest values of $t$, a similar result holds. However, here we find additional extensions or pathological cases, which justify a separated treatment. For $m\geq 4$ consider the Lie algebra $\frak{s}_{m}^{1}$ expressed over the basis $\left\{\omega_{1},..,\omega_{2m-1},\omega_{2m+1},\omega_{2m+2}\right\}$ :
 
\begin{align*}
d\omega_{1}  &  =d\omega_{2}=0\\
d\omega_{j}  &  =\omega_{1}\wedge\omega_{j-1},\;3\leq j\leq2m-1\\
d\omega_{2m+1}  &  =\omega_{2}\wedge\omega_{5}-\omega_{3}\wedge\omega_{4}\\
d\omega_{2m+2}  &  =\omega_{1}\wedge\omega_{2m+1}+2\omega_{2}\wedge\omega
_{6}-\omega_{3}\wedge\omega_{5}%
\end{align*}

\begin{remark}
This algebra plays the same role as $\frak{s}_{m}$ before. Observe also that its characteristic sequence is $\left(2m-2,2,1\right)$. This case corresponds to those models for which the "last" vector is not central.
\end{remark}

\begin{proposition}
Let $\frak{g}$ be a $\left(P2\right)$-algebra of characteristic sequence $\left(2m-1,2,1\right)$ and $h\left(X_{2m+2}\right)=\frac{2t+1}{2}$ with $t=1,2$ over the ordered basis $\left\{X_{1},..,X_{2m+2}\right\}$. Then $\frak{g}$ is a central extension of $\frak{g}_{\left(  m,1\right)  }^{2}$ if $t=1$, a central extension of  $\frak{g}_{\left(  4,2\right)  }^{1}$ if $t=2,m=4$, or a central extension of $\frak{g}_{\left(  m,2\right)  }^{2}$ or $\frak{s}_{m}$ if $t=2$ and $m\geq 4$.
\end{proposition}

\begin{proof}
Again, the main idea of the proof is the same as in theorem 3. We only comment few aspects : for the exceptional ( nine dimensional ) case $\frak{g}_{\left(  4,1\right)  }^{1}$ a central extension satisfying $\left(P2\right)$ is determined by the cocycles  $\varphi_{19}\in H_{4}^{2,\frac{5}{2}}\left(  \frak{g}_{\left(  4,2\right)  }^{3},\mathbb{C}\right)
,\;\varphi_{26},\varphi_{35}\in H_{2}^{2,\frac{5}{2}}\left(  \frak{g}_{\left(
4,2\right)  }^{3},\mathbb{C}\right)  $ and $\varphi_{29}\in H_{5}^{2,\frac
{5}{2}}\left(  \frak{g}_{\left(  4,2\right)  }^{3},\mathbb{C}\right)  $
subjected to the relations
\begin{align*}
\varphi_{26}+2\varphi_{35}  &  =0\\
\varphi_{19}+\varphi_{35}  &  =0\\
2\varphi_{19}+\varphi_{29}  &  =0
\end{align*}
It is clear that they define a unique extension.\newline
Any central extension of degree one of $\frak{s}_{m}^{1}$ is determined by the adjunction of a
differential form, which we will call $d\omega_{2m}.$ The graduation and the
characteristic sequence imply that this differential form is of the type
\[
d\omega_{2m}=\omega_{1}\wedge\omega_{2m-1}+\sum_{i,j}\varphi_{ij}\,\omega
_{i}\wedge\omega_{j}+\varphi_{2,2m+2}\omega_{2}\wedge\omega_{2m+2}%
+\varphi_{3,2m+1}\omega_{3}\wedge\omega_{2m+1}%
\]
where $\varphi_{ij}\in H_{2}^{2,m-1}\left(  G_{m}^{1},\mathbb{C}\right)
,\;\varphi_{2,2m+2},\varphi_{2,2m+1}\in H_{5}^{2,m-1}\left(  G_{m}%
^{1},\mathbb{C}\right)  $, as we have $h\left(  X_{2m}\right)  =m-1$. The
following relations hold
\[
\varphi_{2,2m-1}+\left(  -1\right)  ^{j}\varphi_{j,2m+1-j}=0,\;j=3,..,\left[
\frac{2m+1}{2}\right]
\]%
\[
\varphi_{2,2m+2}+\varphi_{3,2m+1}=0
\]
This implies the existence of a unique extension having characteristic sequence $\left(2m-1,2,1\right)$, and given by the equations
\begin{align*}
d\omega_{1}  &  =d\omega_{2}=0\\
d\omega_{j}  &  =\omega_{1}\wedge\omega_{j-1},\;3\leq j\leq2m-1\\
d\omega_{2m}  &  =\omega_{1}\wedge\omega_{2m-1}+\sum_{j=2}^{\left[
\frac{2m+1}{2}\right]  }\left(  -1\right)  ^{j}\omega_{j}\wedge\omega
_{2m+1-1}+\omega_{2}\wedge\omega_{2m+2}-\omega_{3}\wedge\omega_{2m+1}\\
d\omega_{2m+1}  &  =\omega_{2}\wedge\omega_{5}-\omega_{3}\wedge\omega_{4}\\
d\omega_{2m+2}  &  =\omega_{1}\wedge\omega_{2m+1}+2\omega_{2}\wedge\omega
_{6}-\omega_{3}\wedge\omega_{5}%
\end{align*}
We denote this algebra with $\frak{g}_{\left(  m,2\right)  }^{5}$.
\end{proof}

Now it is not difficult to establish the main result for this characteristic sequence :

\begin{theorem}[Classification of $\left(P2\right)$-algebras of ch.s.
$\left(  2m-1,2,1\right)$]
Let $m\geq 4$ and  $\frak{g}$  be a $\left(2m+2\right)$- dimensional Lie algebra of characteristic sequence $\left(  2m-1,2,1\right)$. Then
$\frak{g}$ is a $\left(P2\right)$-algebra if and only if it is isomorphic to one of the following algebras :

\begin{enumerate}

\item $\frak{g}_{\left(  4,2\right)  }^{1,1}:$%
\begin{align*}
d\omega_{1}  &  =d\omega_{2}=0\\
d\omega_{j}  &  =\omega_{1}\wedge\omega_{j-1},\;3\leq j\leq5\\
d\omega_{6}  &  =\omega_{1}\wedge\omega_{5}+\omega_{2}\wedge\omega_{5}%
-\omega_{3}\wedge\omega_{4}\\
d\omega_{7}  &  =\omega_{1}\wedge\omega_{6}+2\omega_{2}\wedge\omega_{6}%
-\omega_{3}\wedge\omega_{5}-2\omega_{2}\wedge\omega_{9}\\
d\omega_{8}  &  =\omega_{1}\wedge\omega_{7}+\omega_{2}\wedge\omega_{7}%
+\omega_{3}\wedge\omega_{6}-2\omega_{4}\wedge\omega_{5}-2\omega_{3}%
\wedge\omega_{9}\\
d\omega_{9}  &  =\omega_{2}\wedge\omega_{5}-\omega_{3}\wedge\omega_{4}\\
d\omega_{10}  &  =\omega_{1}\wedge\omega_{9}+2\omega_{2}\wedge\omega
_{6}-\omega_{3}\wedge\omega_{5}-2\omega_{2}\wedge\omega_{9}%
\end{align*}

\item $\frak{g}_{\left(  m,0\right)  }^{1+k}\;\left(  k=0,1\right)  :$
\begin{align*}
d\omega_{1}  &  =d\omega_{2}=0\\
d\omega_{j}  &  =\omega_{1}\wedge\omega_{j-1},\;3\leq j\leq2m-1\\
d\omega_{2m}  &  =\omega_{1}\wedge\omega_{2m-1}+\sum_{j=2}^{\left[
\frac{2m+1}{2}\right]  }\left(  -1\right)  ^{j}\,\omega_{j}\wedge
\omega_{2m+1-j}\\
d\omega_{2m+1}  &  =0\\
d\omega_{2m+2}  &  =\omega_{1}\wedge\omega_{2m+1}+k\omega_{2}\wedge
\omega_{2m+1}%
\end{align*}

\item $\frak{g}_{\left(  m,t\right)  }^{2,1}\;\left(  1\leq t\leq m-2\right)
:$%
\begin{align*}
d\omega_{1}  &  =d\omega_{2}=0\\
d\omega_{j}  &  =\omega_{1}\wedge\omega_{j-1},\;3\leq j\leq2m-1\\
d\omega_{2m}  &  =\omega_{1}\wedge\omega_{2m-1}+\sum_{j=2}^{\left[
\frac{2m+1}{2}\right]  }\left(  -1\right)  ^{j}\,\omega_{j}\wedge
\omega_{2m+1-j}\\
d\omega_{2m+1}  &  =\sum_{j=2}^{t+1}\left(  -1\right)  ^{j}\,\omega_{j}%
\wedge\omega_{3-j+2t}\\
d\omega_{2m+2}  &  =\omega_{1}\wedge\omega_{2m+1}+\sum_{j=2}^{t+1}\left(
-1\right)  ^{j}\left(  t+2-j\right)  \,\omega_{j}\wedge\omega_{4-j+2t}%
\end{align*}

\item $\frak{g}_{\left(  m,1\right)  }^{2,2}:$%
\begin{align*}
d\omega_{1}  &  =d\omega_{2}=0\\
d\omega_{j}  &  =\omega_{1}\wedge\omega_{j-1},\;3\leq j\leq2m-1\\
d\omega_{2m}  &  =\omega_{1}\wedge\omega_{2m-1}+\sum_{j=2}^{\left[
\frac{2m+1}{2}\right]  }\left(  -1\right)  ^{j}\,\omega_{j}\wedge
\omega_{2m+1-j}\\
d\omega_{2m+1}  &  =\,\omega_{2}\wedge\omega_{3}\\
d\omega_{2m+2}  &  =\omega_{1}\wedge\omega_{2m+1}+\,\omega_{2}\wedge\omega
_{4}+\omega_{2}\wedge\omega_{2m+1}%
\end{align*}

\item $\frak{g}_{\left(  m,m-2\right)  }^{3,1}:$%
\begin{align*}
d\omega_{1}  &  =d\omega_{2}=0\\
d\omega_{j}  &  =\omega_{1}\wedge\omega_{j-1},\;3\leq j\leq2m-3\\
d\omega_{2m-2}  &  =\omega_{1}\wedge\omega_{2m-3}+\sum_{j=2}^{\left[
\frac{2m+1}{2}\right]  -1}\left(  -1\right)  ^{j}\,\omega_{j}\wedge
\omega_{2m-1-j}\\
d\omega_{2m-1}  &  =\omega_{1}\wedge\omega_{2m-1}+\sum_{j=2}^{m-1}\left(
-1\right)  ^{j}\left(  m-j\right)  \,\omega_{j}\wedge\omega_{2m+1-j}-\left(
m-2\right)  \omega_{2}\wedge\omega_{2m+1}\\
d\omega_{2m}  &  =\omega_{1}\wedge\omega_{2m-1}+\sum_{j=2}^{m}\frac{\left(
-1\right)  ^{j}\left(  j-2\right)  \left(  2m-1-j\right)  }{2}\,\omega
_{j}\wedge\omega_{2m+1-j}\\
d\omega_{2m+1}  &  =\sum_{j=2}^{\left[  \frac{2m+1}{2}\right]  -1}\left(
-1\right)  ^{j}\,\omega_{j}\wedge\omega_{2m-1-j}\\
d\omega_{2m+2}  &  =\omega_{1}\wedge\omega_{2m+1}+\sum_{j=2}^{m-1}\left(
-1\right)  ^{j}\left(  m-j\right)  \,\omega_{j}\wedge\omega_{2m-j}-\left(
m-2\right)  \omega_{2}\wedge\omega_{2m+1}%
\end{align*}

\item $\frak{g}_{\left(  m,2\right)  }^{5}:$%
\begin{align*}
d\omega_{1}  &  =d\omega_{2}=0\\
d\omega_{j}  &  =\omega_{1}\wedge\omega_{j-1},\;3\leq j\leq2m-1\\
d\omega_{2m}  &  =\omega_{1}\wedge\omega_{2m-1}+\sum_{j=2}^{\left[
\frac{2m+1}{2}\right]  }\left(  -1\right)  ^{j}\omega_{j}\wedge\omega
_{2m+1-j}+\omega_{2}\wedge\omega_{2m+2}-\omega_{3}\wedge\omega_{2m+1}\\
d\omega_{2m+1}  &  =\omega_{2}\wedge\omega_{5}-\omega_{3}\wedge\omega_{4}\\
d\omega_{2m+2}  &  =\omega_{1}\wedge\omega_{2m+1}+2\omega_{2}\wedge\omega
_{6}-\omega_{3}\wedge\omega_{5}%
\end{align*}
Moreover, these algebras are pairwise non isomorphic.
\end{enumerate}
\end{theorem}

We proceed stepwise, as done in the previous section.

\begin{proposition}
For $m\geq 4$ the following assertions hold :

\begin{enumerate}
\item $\mathbf{E}_{c,1}\left(  \frak{g}_{\left(  m,1\right)  }^{2},2m-1\right)
=\mathbf{E}_{c,1}^{\frac{3}{2},2,2m-2}\left(  \frak{g}_{\left(  m,1\right)  }^{2},2m-1\right)
+\mathbf{E}_{c,1}^{\frac{3}{2},2,2m-2,2m-1}\left(  \frak{g}_{\left(
m,1\right)  }^{2},2m-1\right)$

\item $\mathbf{E}_{c,1}\left(  \frak{g}_{\left(  m,t\right)  }^{2},2m-1\right)
=\mathbf{E}_{c,1}^{\frac{2t+1}{2},2,2m-2t}\left(  \frak{g}_{\left(  m,1\right)  }^{2},2m-1\right)$

\item $\mathbf{E}_{c,1}\left(  \frak{g}_{\left(  m,m-2\right)  }^{3},2m-1\right)=\mathbf{E}_{c,1}^{\frac{2m-3}{2},2,4,5}\left(  \frak{g}_{\left(  m,m-2\right)  }^{3},2m-1\right)$
\end{enumerate}
\end{proposition}

\begin{proof}
For any case the reasoning is similar to previous ones.

\begin{enumerate}
\item  The cocycle $\varphi_{24}\in H_{2}^{2,\frac{3}{2}}\left(
\frak{g}_{\left(  m,1\right)  }^{2},\mathbb{C}\right)  $ makes reference to
the differential form $d\omega_{2m+1}=\omega_{2}\wedge\omega_{3}$. To this we
have to add, by the characteristic sequence and the closure of the forms
system, the cocycle $\varphi_{1,2m+1}\in H_{2m-2}^{2,\frac{3}{2}}\left(
\frak{g}_{\left(  m,1\right)  }^{2},\mathbb{C}\right)  $, subjected to the
condition $\varphi_{1,2m+1}+\varphi_{24}=0$. A second class of extensions is
defined by the cocycle ( class ) $\varphi_{2,2m+1}\in H_{2m-1}^{2,\frac{3}{2}%
}\left(  \frak{g}_{\left(  m,1\right)  }^{2},\mathbb{C}\right)  $.

\item  The cocycles which define the desired extensions are
\[
\varphi_{j,4-j+2t}\in H_{2}^{2,\frac{2t+1}{2}}\left(  \frak{g}_{\left(
m,t\right)  }^{2},\mathbb{C}\right)  ,\;\varphi_{1,2m+1}\in H_{2m-2t}%
^{2,\frac{2t+1}{2}}\left(  \frak{g}_{\left(  m,t\right)  }^{2},\mathbb{C}%
\right)
\]
satisfying
\begin{align*}
\varphi_{1,2m+1}+\varphi_{t+1,3+t}  &  =0\\
\varphi_{2,2+t}+\left(  -1\right)  ^{j}\left(  t+2-j\right)  \varphi
_{j,4-j+2t}  &  =0,\;3\leq j\leq t+1
\end{align*}

\item  We have to consider the cocycles
\begin{align*}
\varphi_{1,2m+1}  &  \in H_{4}^{2,\frac{2m-3}{2}}\left(  \frak{g}_{\left(
m,m-1\right)  }^{3},\mathbb{C}\right)  ,\;\varphi_{j,2m-j}\in H_{2}%
^{2,\frac{2m-3}{2}}\left(  \frak{g}_{\left(  m,m-1\right)  }^{3}%
,\mathbb{C}\right) \\
\varphi_{2,2m+1}  &  \in H_{5}^{2,\frac{2m-3}{2}}\left(  \frak{g}_{\left(
m,m-1\right)  }^{3},\mathbb{C}\right)
\end{align*}
subjected to the relations
\begin{align*}
\left(  m-2\right)  \varphi_{1,2m+1}+\varphi_{2,2m+1}  &  =0\\
\varphi_{1,2m+1}+\left(  -1\right)  ^{m}\varphi_{m-1,m+1}  &  =0\\
\varphi_{2,2m-2}+\left(  -1\right)  ^{j}\left(  m-j\right)  \varphi_{j,2m-j}
&  =0,\;3\leq j\leq m
\end{align*}
\end{enumerate}
\end{proof}

\begin{remark}
It follows that for other choices of the superindixes, and in particular for fractionary depths, the previous spaces reduce to zero.
\end{remark}

\begin{corollary}
For $m\geq4$ and $1\leq t\leq m-2$ we have
\[
\mathbf{E}_{c,1}^{t,k_{1},..,k_{r}}\left(  \frak{g},2m-1\right)=\{0\}
\]
where $\frak{g}\in\left\{  \frak{g}_{\left(  4,2\right)  }^{1},\frak{g}%
_{\left(  m,m-1\right)  }^{3},\frak{g}_{\left(  m,t\right)  }^{2}\right\}  $.
\end{corollary}

\begin{corollary}
The following identities hold

\begin{enumerate}
\item Any extension $\frak{e}\in\mathbf{E}_{c,1}^{\frac{2m-3}{2},2,4,5}\left(  \frak{g}_{\left(
m,m-2\right)  }^{3},2m-1\right)$ is isomorphic to $\frak{g}_{\left(  m,m-2\right)  }^{3,1}$.

\item Any extension $\frak{e}\in\mathbf{E}_{c,1}^{\frac{3}{2},2,2m-2,2m-1}\left(  \frak{g}_{\left(m,1\right)  }^{2},2m-1\right)$ is isomorphic to  $\frak{g}_{\left(  m,1\right)  }^{2,1}$

\item Any extension $\frak{e}\in\mathbf{E}_{c,1}^{\frac{2t+1}{2},2,2m-2t}\left(  \frak{g}_{\left(m,t\right)  }^{2},2m-1\right)$ is isomorphic to $\frak{g}_{\left(  m,t\right)  }^{2,1}$
\end{enumerate}
\end{corollary}

The proof is elementary.

\bigskip Now we prove the classification theorem :

\begin{proof}[\textbf{Proof of theorem 6} ]
If $h\left(  X_{2m+2}\right)  =\frac{1}{2}$, it is trivial to verify that this
vector must be in the center. Then the factor algebra $\frac{\frak{g}%
}{\left\langle X_{2m+2}\right\rangle }$ has characteristic sequence $\left(
2m-1,1,1\right)  $ and $h\left(  X_{2m+1}\right)  =0$. We
know that for this depth there does not exist any nonsplit model. Thus
$\frac{\frak{g}}{\left\langle X_{2m+2}\right\rangle }\simeq Q_{2m-1}%
\oplus\mathbb{C}$ and $\frak{g}$ \ must be isomorphic to either $\frak{g}%
_{\left(  m,0\right)  }^{1}$ or $\frak{g}_{\left(  m,0\right)  }^{2}$. If
$h\left(  X_{2m+2}\right)  =\frac{3}{2}$ \ the characteristic sequence and the
graduation imply that $X_{2m+2}\in Z\left(  \frak{g}\right)  $, thus
$\frac{\frak{g}}{\left\langle X_{2m+2}\right\rangle }\simeq\frak{g}_{\left(
m,1\right)  }^{2}$. In consequence $\frak{g}\in E_{c,1}\left(  \frak{g}%
_{\left(  m,1\right)  }^{2}\right)  \cap\delta N_{2m}^{2m+2}$.\newline If
$h\left(  X_{2m+2}\right)  =\frac{5}{2}$ $\ $the Jacobi conditions give two
solutions : if $X_{2m+2}\in Z\left(  \frak{g}\right)  $ then $\frac{\frak{g}%
}{\left\langle X_{2m+2}\right\rangle }$ is isomorphic to $\frak{g}_{\left(
m,2\right)  }^{2}$, and if $X_{2m+2}\notin Z\left(  \frak{g}\right)  $ then
$X_{2m}$ must be a central vector, from which \ $\frac{\frak{g}}{\left\langle
X_{2m}\right\rangle }$ is isomorphic to $G_{m}^{1}$; in the first case
$\frak{g}\simeq\frak{g}_{\left(  m,2\right)  }^{2,1}$ and in the second
$\frak{g}\simeq\frak{g}_{\left(  m,2\right)  }^{5}$. \newline For $h\left(
X_{2m+2}\right)  =\frac{2t+1}{2},\;3\leq t\leq m-3$ \ the characteristic
sequence and the graduation imply that $Z\left(  \frak{g}\right)
\supset\left\langle X_{2m+2}\right\rangle $, thus $\frak{g\simeq g}_{\left(
m,t\right)  }^{2}$ by the previous reasoning. Finally, for the depth
$\frac{2m-3}{2}$ the factor of $\frak{g}$ \ by the central ideal $\left\langle
X_{2m}\right\rangle $ is either isomorphic to $\frak{g}_{\left(  m,m-2\right)
}^{2}$ or $\frak{g}_{\left(  m,m-2\right)  }^{3}$.
\end{proof}

The next table resumes the families obtained in theorem $6$ : 
\[
\text{%
\begin{tabular}
[c]{|c|c|c|c|}\hline
\multicolumn{4}{|c|}{Table 1}\\\hline
$\frak{g}$ & $\dim\frak{g}$ & ch.s & type\\\hline
$\frak{g}_{\left(  m,0\right)  }^{1}$ & $2m+2$ & $\left(  2m-1,2,1\right)  $ &
$\left(  3,2,1,..,1\right)  $\\\hline
$\frak{g}_{\left(  m,0\right)  }^{2}$ & $2m+2$ & $\left(  2m-1,2,1\right)  $ &
$\left(  3,2,1,..,1\right)  $\\\hline
$\frak{g}_{\left(  m,t\right)  }^{2,1}$ & $2m+2$ & $\left(  2m-1,2,1\right)  $%
& $\left(  2,1,..,\overset{\left(  2t+1\right)  }{2,}2,1,..1\right)  $\\\hline
$\frak{g}_{\left(  m,1\right)  }^{2,2}$ & $2m+2$ & $\left(  2m-1,2,1\right)  $%
& $\left(  2,1,2,2,1,..,1\right)  $\\\hline
$\frak{g}_{\left(  m,m-2\right)  }^{3,1}$ & $2m+2$ & $\left(  2m-1,2,1\right)
$ & $\left(  2,1,..,2,2,1\right)  $\\\hline
$\frak{g}_{\left(  m,2\right)  }^{5}$ & $2m+2$ & $\left(  2m-1,2,1\right)  $ &
$\left(  2,1,1,1,2,2,1,..,1\right)  $\\\hline
\end{tabular}
}%
\]%

\section{$\left(P2\right)$-algebras of characteristic sequence $\left(2m-1,q,1\right)$}

In this section we describe different families of $\left(P2\right)$ Lie algebras in arbitrary 
dimension and characteristic sequence $\left(2m-1,q,1\right)$ with $q\geq 1$.The algebras we enumerate are obtained by central extensions of the algebras classified in theorem $6$. 
Now observe that for any $q\geq 3$ the classification of $\left(P2\right)$-algebras having the specified characteristic sequence is given up to the exceptional model ( like $\frak{s}_{m}$ and $\frak{s}_{m}^{1}$ before ) which appears for any $q$. The remarkable fact is, however, that for any $q\geq -1$ ( here allowing the cases treated ) most models can be interpreted as central extensions of the algebra $Q_{n}$. This justifies the importance of this model within the $\left(P2\right)$-algebras.

\bigskip Let $m\geq 4$. Consider the Lie algebras 

\begin{itemize}
\item $\frak{g}_{\left(  m,0\right)  }^{1+k,q}\;\left(  k=0,1\right)  ,\;1\leq
q\leq2m-3$%
\begin{align*}
d\omega_{1} &  =d\omega_{2}=0\\
d\omega_{j} &  =\omega_{1}\wedge\omega_{j-1},\;3\leq j\leq2m-1\\
d\omega_{2m} &  =\omega_{1}\wedge\omega_{2m-1}+\sum_{j=2}^{\left[  \frac
{2m+1}{2}\right]  }\left(  -1\right)  ^{j}\omega_{j}\wedge\omega_{2m+1-j}\\
d\omega_{2m+1} &  =0\\
d\omega_{2m+2} &  =\omega_{1}\wedge\omega_{2m+1}+k\omega_{2}\wedge
\omega_{2m+1}\\
d\omega_{2m+2+r} &  =\omega_{1}\wedge\omega_{2m+1+r}+k\omega_{2+r}\wedge
\omega_{2m+1},\;1\leq r\leq q
\end{align*}
over the basis $\left\{\omega_{1},..,\omega_{2m+2+r}\right\}$.

\item $\frak{g}_{\left(  m,t\right)  }^{2,1,q}\;\;\left(  1\leq q\leq
2m-2t-3\right)  $%
\begin{align*}
d\omega_{1} &  =d\omega_{2}=0\\
d\omega_{j} &  =\omega_{1}\wedge\omega_{j-1},\;3\leq j\leq2m-1\\
d\omega_{2m} &  =\omega_{1}\wedge\omega_{2m-1}+\sum_{j=2}^{\left[  \frac
{2m+1}{2}\right]  }\left(  -1\right)  ^{j}\,\omega_{j}\wedge\omega_{2m+1-j}\\
d\omega_{2m+1} &  =\sum_{j=2}^{t+1}\left(  -1\right)  ^{j}\,\omega_{j}%
\wedge\omega_{3-j+2t}\\
d\omega_{2m+2} &  =\omega_{1}\wedge\omega_{2m+1}+\sum_{j=2}^{t+1}\left(
-1\right)  ^{j}\left(  t+2-j\right)  \,\omega_{j}\wedge\omega_{4-j+2t}\\
d\omega_{2m+2+r} &  =\omega_{1}\wedge\omega_{2m+1+r}+\sum_{j=2}^{t+1}\left(
-1\right)  ^{j}\,S_{j}^{r}\,\omega_{j}\wedge\omega_{4-j+2t+r},\;1\leq r\leq q
\end{align*}
where %
\begin{align*}
S_{j}^{1} &  =\sum_{k=j}^{t+1}\left(  t+2-k\right)  ,\;2\leq j\leq t+1\\
S_{j}^{k} &  =\sum_{k=j}^{t+1}S_{j}^{k-1},\;\;2\leq k\leq q
\end{align*}
over the basis $\left\{\omega_{1},..,\omega_{2m+2+r}\right\}$.

\item $\frak{g}_{\left(  m,1\right)  }^{2,2,q}\;\;\left(  1\leq q\leq
2m-3\right)  $%
\begin{align*}
d\omega_{1} &  =d\omega_{2}=0\\
d\omega_{j} &  =\omega_{1}\wedge\omega_{j-1},\;3\leq j\leq2m-1\\
d\omega_{2m} &  =\omega_{1}\wedge\omega_{2m-1}+\sum_{j=2}^{\left[  \frac
{2m+1}{2}\right]  }\left(  -1\right)  ^{j}\,\omega_{j}\wedge\omega_{2m+1-j}\\
d\omega_{2m+1} &  =\,\omega_{2}\wedge\omega_{3}\\
d\omega_{2m+2} &  =\omega_{1}\wedge\omega_{2m+1}+\,\omega_{2}\wedge\omega
_{4}+\omega_{2}\wedge\omega_{2m+1}\\
d\omega_{2m+2+r} &  =\omega_{1}\wedge\omega_{2m+1+r}+\,\omega_{2}\wedge
\omega_{4+r}+\omega_{2}\wedge\omega_{2m+1+r},\;1\leq r\leq q
\end{align*}
over the basis $\left\{\omega_{1},..,\omega_{2m+2+r}\right\}$.

\item $\frak{g}_{\left(  m,2\right)  }^{5,q}\;\left(  1\leq q\leq2m-5\right)
$%
\begin{align*}
d\omega_{1} &  =d\omega_{2}=0\\
d\omega_{j} &  =\omega_{1}\wedge\omega_{j-1},\;3\leq j\leq2m-1\\
d\omega_{2m} &  =\omega_{1}\wedge\omega_{2m-1}+\sum_{j=2}^{\left[  \frac
{2m+1}{2}\right]  }\left(  -1\right)  ^{j}\,\omega_{j}\wedge\omega
_{2m+1-j}+\omega_{2}\wedge\omega_{2m+1}-\omega_{3}\wedge\omega_{2m+1}\\
d\omega_{2m+1} &  =\omega_{2}\wedge\omega_{5}-\omega_{3}\wedge\omega_{4}\\
d\omega_{2m+2} &  =\omega_{1}\wedge\omega_{2m+1}+2\,\omega_{2}\wedge\omega
_{6}-\omega_{3}\wedge\omega_{5}\\
d\omega_{2m+2+r} &  =\omega_{1}\wedge\omega_{2m+1+r}+\left(  2+r\right)
\,\omega_{2}\wedge\omega_{6+r}-\omega_{3}\wedge\omega_{5+r},\;1\leq r\leq q
\end{align*}
over the basis $\left\{\omega_{1},..,\omega_{2m+2+r}\right\}$.

\end{itemize}

\begin{notation}
For $m\geq 4$ and any fixed $q\geq 1$ let $\frak{g}_{q}\in\left\{\frak{g}_{\left(  m,0\right)  }^{1+k,q}, \frak{g}_{\left(  m,t\right)  }^{2,1,q}, \frak{g}_{\left(  m,1\right)  }^{2,2,q}, \frak{g}_{\left(  m,2\right)  }^{5,q}\right\}$.
\end{notation}

\begin{theorem}
For $q\geq 1$ the Lie algebra $\frak{g}_{q}$ is a central extension of $\frak{g}_{q-1}$ by $\mathbb{C}$. Moreover, $\frak{g}_{q}$ is a $\left(P2\right)$-algebra of characteristic sequence  
$\left(2m-1,2+q,1\right)$.
\end{theorem}

\begin{proof}
We prove the assertion for $\frak{g}_{\left(  m,t\right)  }^{2,1,q}$. For the remaining cases the reasoning is similar.\newline
Recall that for $\frak{g}_{\left(m,t\right)}^{2,1}$ the last differential form is given by
\[
d\omega_{2m+2}   =\omega_{1}\wedge\omega_{2m+1}+\sum_{j=2}^{t+1}\left(
-1\right)  ^{j}\left(  t+2-j\right)  \,\omega_{j}\wedge\omega_{4-j+2t}
\]
A central extension of $\frak{g}_{\left(m,t\right)}^{2,1}$ by $\mathbb{C}$ which is a $\left(P2\right)$-algebra will be determined by the adjunction of a differential form $d\omega_{2m+3}$, whose structure is 
\[
d\omega_{2m+3}   =\omega_{1}\wedge\omega_{2m+2}+\sum_{j=2}^{t+1}\left(
-1\right)  ^{j}\,\varphi_{j,5-j+2t}\,\omega_{j}\wedge\omega_{5-j+2t},
\]
where the cocycles
\begin{equation*}
\varphi _{j,5-j+2t}\in H_{2}^{2,t+1}\left( \frak{g}_{\left( m,t\right)
}^{2,1},\mathbb{C}\right) 
\end{equation*}
satisfy
\begin{equation*}
\varphi _{2,3+2t}+\left( -1\right) ^{j}\sum_{k=j}^{t+1}\left( t+2-k\right)
\varphi _{j,5-j+2t}=0,\;3\leq j\leq t+1
\end{equation*}
We thus obtain a unique extension class which is isomorphic to $\frak{g}%
_{\left( m,t\right) }^{2,1,1}$. This shows the assertion for $q=1$. Let it
be true for $q_{0}>1$. Then the Cartan-Maurer equations of $\frak{g}_{\left(
m,t\right) }^{2,1,q_{0}}$ are 
\begin{equation*}
\begin{align*}
d\omega_{1} &  =d\omega_{2}=0\\ 
d\omega_{j} &  =\omega_{1}\wedge\omega_{j-1},\;3\leq j\leq2m-1\\ 
d\omega_{2m} &  =\omega_{1}\wedge\omega_{2m-1}+\sum_{j=2}^{\left[  \frac 
{2m+1}{2}\right]  }\left(  -1\right)  ^{j}\,\omega_{j}\wedge\omega_{2m+1-j}\\ 
d\omega_{2m+1} &  =\sum_{j=2}^{t+1}\left(  -1\right)  ^{j}\,\omega_{j}% 
\wedge\omega_{3-j+2t}\\ 
d\omega_{2m+2} &  =\omega_{1}\wedge\omega_{2m+1}+\sum_{j=2}^{t+1}\left( 
-1\right)  ^{j}\left(  t+2-j\right)  \,\omega_{j}\wedge\omega_{4-j+2t}\\ 
d\omega_{2m+2+r} &  =\omega_{1}\wedge\omega_{2m+1+r}+\sum_{j=2}^{t+1}\left( 
-1\right)  ^{j}\,S_{j}^{r}\,\omega_{j}\wedge\omega_{4-j+2t+r},\;1\leq r\leq q_{0}
\end{align*}
\end{equation*} 
where % 
\begin{equation*}
\begin{align*} 
S_{j}^{1} &  =\sum_{k=j}^{t+1}\left(  t+2-k\right)  ,\;2\leq j\leq t+1\\ 
S_{j}^{k} &  =\sum_{k=j}^{t+1}S_{j}^{k-1},\;\;2\leq k\leq q 
\end{align*} 
\end{equation*}
Now we extend this algebra by $\mathbb{C}$. Supposing that the extension
satisfies the centralizer property and is naturally graded of the prescribed
characteristic sequence, the determining cocycles are 
\begin{eqnarray*}
\varphi _{j,4-j+2t+q_{0}+1} &\in &H_{2}^{2,\frac{2t+2+r}{2}}\left( \frak{g}%
_{\left( m,t\right) }^{2,1,q_{0}},\mathbb{C}\right) \;\text{if\ }r\equiv
1\left( mod\;2\right)  \\
\varphi _{j,4-j+2t+q_{0}+1} &\in &H_{2}^{2,t+\frac{r}{2}+1}\left( \frak{g}%
_{\left( m,t\right) }^{2,1,q_{0}},\mathbb{C}\right) \;\text{if\ }r\equiv
0\left( mod\;2\right) 
\end{eqnarray*}
We have the relations 
\begin{equation*}
\varphi _{j,4-j+2t+q_{0}+1}+\left( -1\right)
^{j}\sum_{j=2}^{t+1}S_{j}^{q_{0}}\omega _{j}\wedge \omega
_{5-j+2t+q_{0}}=0,\;3\leq j\leq t+1
\end{equation*}
and by an elementary change of basis, the adjoined differential form $%
d\omega _{2m+3+q_{0}\text{ }}$ is of type 
\begin{equation*}
d\omega _{2m+3+q_{0}}=\omega _{1}\wedge \omega
_{2m+2+q_{0}}+\sum_{j=2}^{t+1}\left( -1\right) ^{j}S_{j}^{q_{0}}\omega
_{j}\wedge \omega _{5-j+2t+q_{0}}
\end{equation*}
Both the characteristic sequence and centralizer property are obviously satisfied.
\end{proof}

The algebras $\frak{g}_{\left(  m,m-2\right)  }^{3,1}$ only admit one more extension which is a $\left(P2\right)$-algebra. This is due to the extremal position of the vectors that give the two dimensional Jordan block of the characteristic sequence.

\begin{proposition}
For $m\geq 4$ the algebra $\frak{g}_{\left(  m,m-2\right)  }^{3,1,1}$ given by 
\begin{align*}
d\omega_{1} &  =d\omega_{2}=0\\
d\omega_{j} &  =\omega_{1}\wedge\omega_{j-1},\;3\leq j\leq2m-3\\
d\omega_{2m-2} &  =\omega_{1}\wedge\omega_{2m-3}+\sum_{j=2}^{\left[
\frac{2m+1}{2}\right]  -1}\left(  -1\right)  ^{j}\,\omega_{j}\wedge
\omega_{2m-1-j}\\
d\omega_{2m-1} &  =\omega_{1}\wedge\omega_{2m-1}+\sum_{j=2}^{m-1}\left(
-1\right)  ^{j}\left(  m-j\right)  \,\omega_{j}\wedge\omega_{2m+1-j}-\left(
m-2\right)  \omega_{2}\wedge\omega_{2m+1}\\
d\omega_{2m} &  =\omega_{1}\wedge\omega_{2m-1}+\sum_{j=2}^{m}\frac{\left(
-1\right)  ^{j}\left(  j-2\right)  \left(  2m-1-j\right)  }{2}\,\omega
_{j}\wedge\omega_{2m+1-j}-\left(  m-2\right)  \omega_{3}\wedge\omega_{2m+1}\\
d\omega_{2m+1} &  =\sum_{j=2}^{\left[  \frac{2m+1}{2}\right]  -1}\left(
-1\right)  ^{j}\,\omega_{j}\wedge\omega_{2m-1-j}\\
d\omega_{2m+2} &  =\omega_{1}\wedge\omega_{2m+1}+\sum_{j=2}^{m-1}\left(
-1\right)  ^{j}\left(  m-j\right)  \,\omega_{j}\wedge\omega_{2m-j}-\left(
m-2\right)  \omega_{2}\wedge\omega_{2m+1}\\
d\omega_{2m+3} &  =\omega_{1}\wedge\omega_{2m+2}+\sum_{j=2}^{m-1}\left(
-1\right)  ^{j}S^{j}\omega_{j}\wedge\omega_{2m+1-j}-\left(  m-2\right)
\omega_{3}\wedge\omega_{2m+1}%
\end{align*}
over the basis $\left\{\omega_{1},..,\omega_{2m+3}\right\}$ is a $\left(P2\right)$-algebra of characteristic sequence $\left(2m-1,3,1\right)$.
\end{proposition}

\begin{proof}
Any extensions which satisfies the centralizer property, preserves the graduation and has characteristic sequence $\left(2m-1,3,1\right)$ is determined by

\begin{align*}
\varphi_{1,2m+2}  &  \in H_{4}^{2,m-1}\left(  \frak{g}_{\left(
m,m-1\right)  }^{3,1},\mathbb{C}\right)  ,\;\varphi_{j,2m+1-j}\in H_{2}%
^{2,m-1}\left(  \frak{g}_{\left(  m,m-1\right)  }^{3,1}%
,\mathbb{C}\right) \\
\varphi_{3,2m+1}  &  \in H_{5}^{2,m-1}\left(  \frak{g}_{\left(
m,m-1\right)  }^{3,1},\mathbb{C}\right)
\end{align*}
subjected to the relations
\begin{align*}
\left(  m-2\right)  \varphi_{1,2m+2}+\varphi_{3,2m+1}  &  =0\\
\varphi_{1,2m+2}+\left(  -1\right)  ^{m}\varphi_{m-1,m+2}  &  =0\\
\varphi_{2,2m-2}+\left(  -1\right)  ^{j}S^{j} \varphi_{j,2m+1-j}
&  =0,\;3\leq j\leq m
\end{align*}
where $S^{j}=\sum_{j=2}^{m-1}\left(m-j\right)$.\newline Then the class is unique, and by an elementary change of basis the extended algebra is easily seen to be isomorphic to $\frak{g}_{\left(  m,m-2\right)  }^{3,1,1}$. The centralizer property is given by the form $d\omega_{2m}$.
\end{proof}

We resume the result in the following table :
\[%
\begin{tabular}
[c]{|c|c|c|c|}\hline
\multicolumn{4}{|c|}{Table 2}\\\hline
$\frak{g}$ & $\dim\frak{g}$ & ch.s & type\\\hline
$\frak{g}_{\left(  m,0\right)  }^{1,q}$ & $2m+2+q$ & $\left(
2m-1,2+q,1\right)  $ & $\left(  3,2,..,\overset{\left(  2+q\right)  }%
{2},1,..,1\right)  $\\\hline
$\frak{g}_{\left(  m,0\right)  }^{2,q}$ & $2m+2+q$ & $\left(
2m-1,2+q,1\right)  $ & $\left(  3,2,..,\overset{\left(  2+q\right)  }%
{2},1,..,1\right)  $\\\hline
$\frak{g}_{\left(  m,t\right)  }^{2,1,q}$ & $2m+2+q$ & $\left(
2m-1,2+q,1\right)  $ & $\left(  2,1,..,\overset{\left(  2t+1\right)  }%
{2},2,..,\overset{\left(  q+2t+1\right)  }{2},1,..,1\right)  $\\\hline
$\frak{g}_{\left(  m,1\right)  }^{2,2,q}$ & $2m+2+q$ & $\left(
2m-1,2+q,1\right)  $ & $\left(  2,1,2,..,\overset{\left(  3+q\right)  }%
{2},1,..,1\right)  $\\\hline
$\frak{g}_{\left(  m,2\right)  }^{5,q}$ & $2m+2+q$ & $\left(
2m-1,2+q,1\right)  $ & $\left(  2,1,1,1,2,..,\overset{\left(  5+q\right)  }%
{2},1,..,1\right)  $\\\hline
$\frak{g}_{\left(  m,m-2\right)  }^{3,1,1}$ & $2m+3$ & $\left(
2m-1,3,1\right)  $ & $\left(  2,1,..,1,2,2,2\right)  $\\\hline
\end{tabular}
\]

\begin{remark}
Finally, the pathological case $\frak{g}_{\left(4,2\right)}^{1,1}$ admits the extension
$\frak{g}_{\left(  4,2\right)  }^{1,1,1}$ given by 

\begin{align*}
d\omega_{1}  & =d\omega_{2}=0\\
d\omega_{j}  & =\omega_{1}\wedge\omega_{j-1};\;3\leq j\leq5\\
d\omega_{6}  & =\omega_{1}\wedge\omega_{5}+\omega_{2}\wedge\omega_{5}%
-\omega_{3}\wedge\omega_{4}\\
d\omega_{7}  & =\omega_{1}\wedge\omega_{6}+2\omega_{2}\wedge\omega_{6}%
-\omega_{3}\wedge\omega_{5}-2\omega_{2}\wedge\omega_{9}\\
d\omega_{8}  & =\omega_{1}\wedge\omega_{7}+\omega_{2}\wedge\omega_{7}%
-\omega_{3}\wedge\omega_{6}+2\omega_{4}\wedge\omega_{5}-2\omega_{3}%
\wedge\omega_{9}\\
d\omega_{9}  & =\omega_{2}\wedge\omega_{5}-\omega_{3}\wedge\omega_{4}\\
d\omega_{10}  & =\omega_{1}\wedge\omega_{9}+2\omega_{2}\wedge\omega_{6}%
-\omega_{3}\wedge\omega_{5}-2\omega_{2}\wedge\omega_{9}\\
d\omega_{11}  & =\omega_{1}\wedge\omega_{10}+3\omega_{2}\wedge\omega
_{7}-\omega_{3}\wedge\omega_{6}-2\omega_{3}\wedge\omega_{9}%
\end{align*}
As both the dimension and the characteristic sequence [$\left(7,3,1\right)$] are fixed, this algebra is not of great interest for the general case.
\end{remark}

\textit{ }
\end{document}